\documentclass[12pt,amssymb]{article}
\usepackage{amsmath}
\usepackage{amsthm}
\usepackage{amssymb}
\usepackage{graphicx}
\newcommand{\DETAILS}[1]{}
\newtheorem{remark}{Remark}

\begin{document}

\newcommand{\DATUM}{September 15, 2009}              
\pagestyle{myheadings}                         
\markboth{\hfill{WP, \DATUM}}{{WM, \DATUM}\hfill}  %

\title{On Collapse of Wave Maps}
\author{
Yu.~N.~Ovchinnikov\
\thanks{ 
Supported under the grant 
RFBR-07-0-2-12058}\\
\small{Max-Planck Institute for Physics of Complex Systems, 01187 Dresden, Germany} \\
\small{L.D.~Landau Inst.~for Theoretical Physics, 
Chernogolovka, 142432, Russia} \\
I. M. Sigal\
\thanks{Supported by NSERC Grant NA7601} \\
\small{Department~of Mathematics, University of Toronto, Toronto, Canada}
 }

\date{\DATUM}

\maketitle

\begin{abstract}
  We derive the universal collapse law of degree 1
equivariant wave maps (solutions of the sigma-model) from the 2+1 Minkowski space-time,
 to the 2-sphere. To this end we introduce a nonlinear transformation from original variables to blowup ones. Our formal derivations are confirmed by numerical simulations.

\end{abstract}

\bigskip




\section{ Introduction}

In this paper we investigate the phenomenon of collapse of degree $1$
equivariant wave maps (solutions of the $\sigma$-model) from the $2+1$ Minkowski space-time,
${\mathbb{R}}^{2+1}$, to the $2$-sphere, $S^2$.  Besides of purely
mathematical interest, the study of the blowup phenomena for such
maps is motivated by the recent efforts to understand the
singularity formation in general relativity \cite{R}.

A wave map, $\varphi$, from a $(d+1)-$dimensional space-time, $M$, with a
metric $\eta$ to a Riemannian manifold, $N$, with a metric $g$, is a
critical point of the action functional, which, in local coordinates, has the form
$$S(\varphi):=\frac{1}{2}\int g_{AB}\partial_{a}\varphi^A \partial_{b}\varphi^B
\eta^{ab}\sqrt{-\eta}d^{d+1}x.$$
(This action functional is also known as the $\sigma-$ model). Critical points of $S(\varphi)$ satisfy the Euler-lagrange equation
$$\partial_{a}\partial^{a}\varphi^A +\Gamma_{BC}^{A}(\varphi)\partial_{a}\varphi^B\partial^{a}\varphi^C=0,$$
where $\Gamma_{BC}^{A}(\varphi)$ is the Christoffel symbols on  $N$. This system of nonlinear PDEs is Hamiltonian, and in particular has conserved energy, $\mathcal{E}(\varphi)$, and scale invariant in the sense that if $\varphi(x)$ is a solution then so is $\varphi(\lambda x)$. The energy, $\mathcal{E}(\varphi)$,  is transformed under scaling as
$$\mathcal{E}(\varphi_{\lambda})= \lambda^{2-d}\mathcal{E}(\varphi),$$ where $\varphi_{\lambda}(x)= \varphi(\lambda x)$. Thus the case $d=2$ of interest for us is the energy critical.

For a map, $\varphi$, to have finite energy, it should converge to a constant
at infinity. In this case for each moment of time, $t$, $\varphi$
can be extended to a continuous map from $S^d$ to $N$ taking the
point at infinity to the limit of $\varphi (x)$ at the spatial infinity. Then one can define the
degree, $\textrm{deg} \varphi$, 
as the homotopy class of $\varphi$ as a map
from $S^d$ to $N$. This degree is conserved under the dynamics generated
by the Euler-Lagrange equations above.

In our case, $M$ is the $2+1$ Minkowski space-time,
${\mathbb{R}}^{2+1}$ and $N$ is the $2$-sphere, $S^2$ with the
standard metric $g:=du^2 + \sin^2u d\theta^2$. In this case the
degree of $\varphi$ is an integer (the degree for maps from $S^2$ to
$S^2$).  Moreover, one has the Bogomolnyi inequality (\cite{B})
$$\mathcal{E}(\varphi) \ge 8\pi |\textrm{deg} \varphi|,$$
which leads to self-dual and anti-dual equations for the minimizers of the static energy for fixed degrees. These equations have explicit solutions (harmonic or anti-harmonic maps). These solutions will be written out below.

Among the maps of the degree $k$ the simplest, most
symmetric maps are the 'radially symmetric' or equivariant maps which are of the
form $\varphi_k(\rho,\phi,t)= (u_k(\rho,t), k\phi)$, where $(\rho,\phi)$ are
the polar coordinates  in ${\mathbb{R}}^2$ and  $(\varphi,\theta)$ are
the spherical coordinates in $S^2$. 
Then the Euler-Lagrange equation for $\varphi$
reduces to
the equation:
\begin{equation}
\label{1} \ddot u = \Delta u - \frac{k^2}{2 \rho^2} \sin (2u)
\end{equation}
for $u_k$. Here $\Delta$ is the $2D$ spherical Laplacian. Moreover,
$\textrm{deg} \varphi_k = Q(u_k) = k$, where $$Q(u):=
\frac{1}{\pi}(u(\infty) - u(0)).$$

Numerical studies of Eqn \eqref{1} led to a conjecture that large-energy, degree one
initial data develop singularities in finite time and the
singularity formation has the universal form of adiabatic shrinking
of the degree-one harmonic map from $\mathbb{R}^2$ to $S^2$  \cite{BCT}. Later, it
was shown by Struwe \cite{Str1} that the existence of a nontrivial harmonic
map is in fact the necessary condition for blowup for $2+1$
equivariant wave maps. In this paper we address the
question of the dynamics of the blow-up process. We show that there is $0 < t_* < \infty$ such that, as $t \rightarrow t_*$, we have on bounded domains in  $\mathbb{R}^2$
$$u(\rho, t) \approx U(\rho/\lambda(t)),$$
where $U(\rho)$ is the profile of the degree $1$ equivariant, static (and in particular harmonic) map, minimizing static energy (see below), and the scaling parameter $\lambda (t)$, satisfies the following second order ODE:
\begin{equation}
\label{16}
 \lambda\ddot\lambda = \frac{\dot\lambda^2}
{\ln(\frac{a}{\lambda\ddot \lambda })},\ \mbox{with}\ a= (1.04)^2 e^{-2}\approx .146.
\end{equation}
We expect that, proceeding as in \cite{BOS}, we can show that the error term in the above relation is $O(\dot\lambda^2)$.

Note that Eqn \eqref{16} shows that if
$\dot\lambda|_{t=0}<0$, then $\dot\lambda<0,\ \ddot\lambda>0$ for $t>0$ and therefore $\lambda$ 
and $|\dot\lambda|$ decrease as $t \to t_*$. Since $\dot\lambda^2$ is the small parameter in our analysis (adiabatic regime),  our
approximation improves as $t\to t_*$.

 An approximate solution of Eqn  \eqref{16} with two free parameters (constants of integration), $t^*$ and $c$, is (see Section \ref{eqn2} below) 
\begin{equation}
\label{17}
 \sqrt{a}(t^{*} - t) = \lambda e^  {\ln^{1/2}
 ( \frac{c}{\lambda}  ) } +
c\frac{\sqrt{\pi}}{2}  e^{1/4} \biggl [ 1 - \Phi  ( -1/2 + \ln^{1/2}
 ( \frac{c}{\lambda}  )  ) \biggr ],
\end{equation}
where  $\Phi(x)\equiv \mbox {\rm erf}(x)$  is the Fresnel integral \cite{Handbook}.
An exact solution of Eqn  \eqref{16} is obtained in Section \ref{eqn2} (see Eqn  \eqref{81}).  A
comparison of the leading term  of this solution with a numerical solution of Eqn \eqref{1} is
given in Fig 1.  This figure shows that the two resulting curves are indistinguishable for times
sufficiently close to the blow-up time.

\begin{figure}[h]
\includegraphics[width=\textwidth]{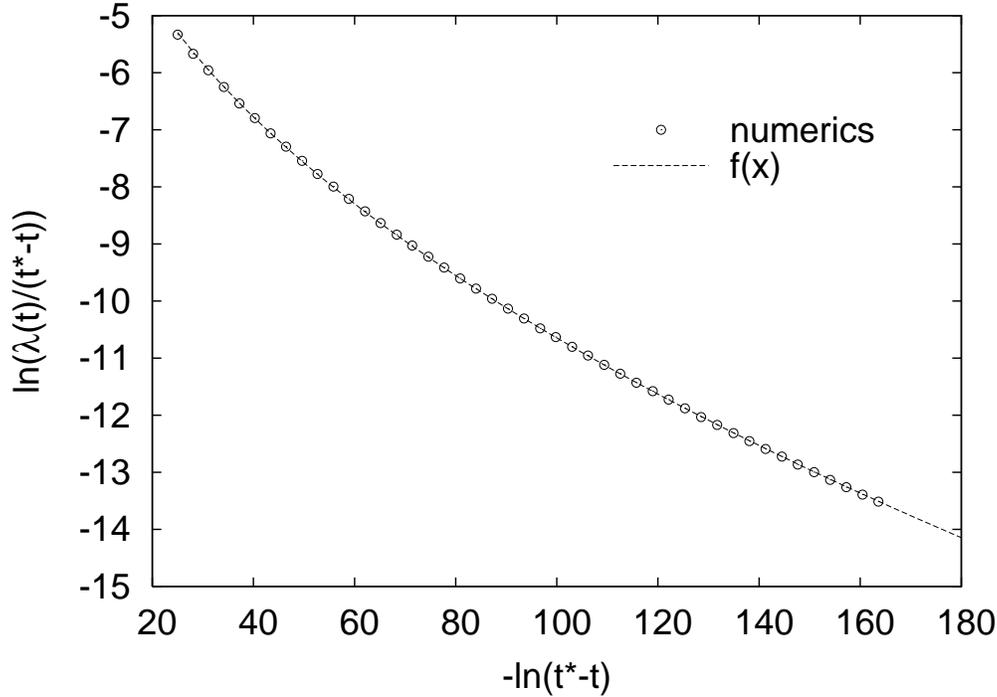}
\caption{\small{For a numerical solution that blows up at time $t^*$ we plot
$y=\ln{\dfrac{\lambda(t)}{t^*-t}}$ as a function of $x=-\ln(t^*-t)$ (circles) and compare it with
the analytic formula $y=f(x)=\dfrac{1}{2}\ln(a)-\sqrt{x+b}$, where $a=0.146$ and $b$ is a free
(non-universal) parameter. Fitting $b$ we get an excellent agreement between numerical and
analytical results.}} \label{figl}
\end{figure}

Observe that 
like Eqn
\eqref{1}, Eqn \eqref{16} is a Hamiltonian equation. Its Lagrangian is
\begin{equation}
\label{19}
 L:= h (\dot\lambda) - \ln \lambda,
\end{equation}
 where the function $h$ is defined by $h'' (x) = -1/
{f^{-1} (x^2/a)}$ with $f (x) = x \ln(1/x) $ (see Section \ref{Ham}).

 The local well-posedness for the wave map equations in Sobolev spaces was proven in \cite{KlMach,KlSel1, KlSel2}, while the global well-posedness for small initial conditions, in \cite{Kov, Sid, Tat1, Tat2, Tao1, Tao2} (see also \cite{CS, CTZ1, CTZ2, Kr1, Kr2, Sh, ShTZ1, Sid, Tat3, Tat4}).
The research on the problem of blowup for the wave maps started with numerical work  \cite{BCT, LPZ,  PZ}. (We do not review here related works for nonlinear wave equations.)

The first numerical evidence for singularity formation for $2+1$
equivariant wave maps to the 2-sphere was given in \cite{BCT}.
In this paper (concerned only with $k=1$ homotopy) the authors showed that blowup
has the form of adiabatic shrinking of the harmonic map and
formulated conjectures about blowup for large energy, blowup
profile and energy concentration and  
that $\lambda(t)/(T-t)$ must go to zero.
 As was already mentioned, it was shown rigorously in \cite{Str1} that the existence of a stationary solution is a necessary condition for the blowup to take place. The blowup scenarios were further numerically investigated in \cite{IL,  LS} (see references therein for additional works).

The first rigorous results on the blowup rate and profile were obtained in   \cite{RS, KST2}. In particular,  \cite{RS}
has obtained the lower bound on the contraction rate $\lambda(t)$ for $k\ge 4$ wave maps. As it turned out this lower bound conforms exactly to the dynamical law derived for the $4+1$ Yang-Mills $k=1$ equivariant solutions in \cite{BOS}  using a formal but careful analysis, explained below in this introduction,  justified by numerical computations. (Earlier numerical analysis for the latter model was
announced in \cite{BT} and described more completely  in the survey \cite{Bi}.) (It was noticed in \cite{RS} (see below), that the $k\ge 2$ wave map equation are similar to the $4+1$ Yang-Mills one for $k=1$. )
Finally, for each $b>1/2$, \cite{KST2} has constructed special solutions of the $k=1$ equivariant wave map equation, Eqn \eqref{1} with $k=1$, which blow up at the rate $\lambda(t) \sim (T-t)^{1+b}$.

Equation \eqref{1} belongs to a general class of semilinear wave
equations in ${\mathbb{R}}^{2+1}$ of the form
\begin{equation}
\label{2} \ddot u = \Delta u - \frac{1}{\rho^2} f(u).
\end{equation}
In the case of 
$f(u) = \frac {k^2}{ 2} \sin (2u)$, Eqn \eqref{2} is, as
was already mentioned, the equation for the profile of the
equivariant wave map from the $2+1$ Minkowski space-time of degree $k$,
${\mathbb{R}}^{2+1}$, to the $2$-sphere, $S^2$. More generally, \eqref{2} is
satisfied by equivariant maps for the case when $N$ is the surface
of revolution with the metric $g:=du^2 + g^2(u) d\theta^2$, where
$g(u)$ is related to $f(u)$ as $f(u)=g(u)g'(u)$.

In the case of $f(u) = 2u (u^2 -1)$
the corresponding equation,
\begin{equation}
\label{3}  \ddot u = \Delta u - \frac{2}{\rho^2}  (u^2 -1)u,
\end{equation} is related to the equation for
equivariant Yang-Mills fields of degree $1$ in the $4+1$ dimensions.

Note that

 (i)
  Eqn \eqref{2} is invariant
 with respect to the scaling transformation, $$u(\rho, t)\rightarrow
 u(\frac{\rho}{\lambda}, \frac{t}{\lambda});$$

 (ii)
  Eqn \eqref{2} can be presented as a Hamiltonian system with the standard symplectic form
 and the Hamiltonian
 \begin{equation}
 \label{2_1}
 H(u,v) :=\int_0^\infty (\frac{1}{2} v^2 + \frac{1}{2}|\nabla u|^2 +  \frac{1}{\rho^2}F(u))\rho d\rho,
  \end{equation}
with $F'(u)=f(u)$.
The scaling properties of the Hamiltonian $H(u,v)$ imply that the
dimension $d=2$ is the critical dimension for Eqn \eqref{2}. This is the
dimension treated in this paper.

We assume now that $f(u)$ is a derivative of a
\textit{double-well potential} $F(u)$, i.e. $F(u)$ is nonnegative and has at
least two global minima, say at $a$ and $b$ for some $b>a$, with
$F(a) =F(b)= 0$, and no minima between $a$ and $b$
($F(u)=\frac{1}{2}\sin^2u$ and $a= 0$, $b= \pi$ in the case of $f(u)
= \frac {k^2}{2} \sin (2u)$ and $F(u)=\frac{1}{2}(u^2-1)^2$ and $a= -1$,
$b=1$ in the case of $f(u) = 2u (u^2-1)$). In this case Eqn \eqref{2} has
the following features:
\begin{enumerate}
 \item[(A)] For each $k \in \mathbb{N}$, Eqn \eqref{2} has static solutions, $U_k(\rho)$ and
$U_{-k}(\rho)=-U_k(\rho)$; they have topological degrees $Q(U_k) = k$ and $Q(U_{-k}) = - k$;
 \item[(B)] For $k =1$ , the solution $U_1(\rho)$ is monotonically increasing from $a$ to
 $b$,  while $U_{-1}(\rho)$ is monotonically
decreasing from $a$ to $b$;
 \item[(C)] The solution, $U_k(\rho)$, is a minimizer of the static energy
functional $E(u)$ under the constrain, $Q(u)=k$, on the topological
charge;
 \item[(D)]  Eqn \eqref{2} conserves the topological charge $Q(u):=  \frac{1}{b-a}(u(\infty) - u(0))$.
\end{enumerate}

Existence of the solutions $U_k(\rho)$ follows from the Bogomolnyi argument, see above. The solutions $U_1(\rho)$ and $U_{-1}(\rho)$ are  called the kink
solution and anti-kink solution, or simply \textit{kink} and
\textit{anti-kink}, respectively.  Since the analysis for $k<0$ can
be obtained from analysis for the case $k>0$ by simply flipping the
signs, in what follows we assume that $k>0$.
Note that though Eqn \eqref{2} is scale
invariant, its static kink solution $U_k(\rho)$ are not. Hence Eqn \eqref{2}
has an entire family, $U_k(\frac{\rho}{\lambda})$, of kink solutions
(symmetry breaking).

From now on \textit{we concentrate on the kink solution}, $U_1(\rho)$ and omit the subindex
$1$: $U_1(\rho) \equiv U(\rho)$.

There is a feature of Eqn \eqref{2} which is not apparent at the first
sight but which plays an important role in our analysis of the
collapse.  The fact that the kink, $U(\rho)$, breaks scale
invariance manifests itself in appearance of the dilation zero
mode
$$\zeta(\rho) := \frac{1}{2}\rho \partial_{\rho}U(\rho).$$
This is a zero eigenfunction, $L_{\rho}\zeta= 0$, for the
linearization operator
\begin{equation}
\label{4}  L_{\rho}=- \frac 1\rho ~ \frac{\partial}{\partial \rho}
\biggl ( \rho~\frac{\partial}{\partial \rho} \biggr ) +
\frac{1}{\rho^2}f'(U(\rho))
\end{equation}
(negative Fr\'echet derivative) for the r.h.s. of \eqref{2} around
$U(\rho)$. This zero mode presents an obstruction to solving Eqn
(2) perturbatively, starting with $U(\rho)$, which can be resolved
by a modulation theory, provided $\zeta$ is an $L^2$ function,
i.e. one can use a Hilbert space spectral theory.

Thus equations of the form \eqref{2} can be organized in two classes
according to which of the following two properties takes place
\begin{enumerate}
 \item[(i)] $\zeta\ is\ in\ L^2$
 \item[(ii)] $\zeta\ is\ not\ in\ L^2$.
\end{enumerate}
The Yang-Mills equation,  \eqref{3}, belongs to the first class while
the wave map equation, \eqref{1}, belong to the second. Indeed, the kink
solution for  \eqref{3} is $U(\rho)=\frac{1-\rho^2}{1+\rho^2}$ and
the corresponding zero mode is $\zeta (\rho)=
\frac{4\rho^2}{(1+\rho^2)^2}$ (see \cite{ Bi, BOS}). For Eqn \eqref{1} with $k=1$ the kink
solution is
$$
U(\rho) = 2~\mbox {\rm  arctan}~ \rho,
$$
while the scaling zero mode is
$$
\zeta (\rho) := \frac{1}{2}\rho\partial_{\rho} U(\rho) =
\frac{\rho}{1+\rho^2}.
$$
Clearly, $\zeta$ is  $L^2$ in the former case and is not $L^2$ in the latter case. (This is possible
due to the fact that 
the operator $L$ has no spectral gap: $\sigma(L)=\sigma_{cont}(L)=[0,\infty)$ (see below). The fact that there is a problem with the modulation approach
due to the nonintegrability of  the zero mode was pointed out by P.Bizo\'n in 2001, \cite{BiPC}.)

Note that $U_k(\rho)=2 \arctan(\rho^k)$
and the corresponding zero mode is square-integrable for $k>1$. Thus in this case we expect that at least the formal analysis of \cite{BOS} of the Yang-Mills equation should go through. (Higher degree equivariant, static solutions for the Yang-Mills equations in $4+1$ dimensions, known as instantons, can be found in \cite{At, Schw}.)

We are interested in solution with initial conditions near the kink
manifold $$M_{kink}:= \{U(\rho/\lambda)| \lambda >0\}.$$ The
fact that the square integrability of of the zero mode $\zeta$ plays
an important role in analysis of such solutions can be gleaned from
the observation that the effective action $S(U_\lambda)$ on the
family $U_\lambda (\rho) := U (\rho/\lambda)$ (the 'effective
action' of $\lambda$) is equal to $$S(U_\lambda)= \frac{1}{2}\int
\{(\lambda\dot{ \lambda})^2 \|\zeta\|^2 - V(U)\} dt,$$ where $V(u):= \frac{1}{2}|\nabla u|^2 +  \frac{1}{\rho^2}F(u),$ with $F'(u)=f(u)$, and diverges,
if $\zeta \not \in L^2$. (For a connection to the geodesic hypothesis see \cite{MS, W}). Here and in what follows $\dot{ \lambda} =\frac{\partial\lambda}{\partial
t}$.

We present heuristic arguments motivating our approach. It is natural to guess that for an initial condition close to the
manifold $M_{kink}$ the solution evolves along this manifold. Let
$U(\frac{\rho}{\lambda(t)})$ be the projection of the solution on
this manifold. If for this projection $\lambda(t)\rightarrow 0$ as
$t \rightarrow t_*$ for some $t_*$, then the solution collapses at
the time $t_*$. With this in mind we look for solutions to Eqn \eqref{2}
in the form
\begin{equation}
\label{5} u(\rho,t) =U(x) + w(x,t),
\end{equation}
where $x=\rho/\lambda$,  a blow-up variable, with $\lambda$ a slowly varying function of time $t$
(we do not pass to the blow-up time variable). Note that while in a
standard approach the scaling, $\lambda$, is fixed at the very beginning
(with corrections at certain scales possibly considered later on) we
leave it free and we look for a differential equation for $\lambda$
which guarantees that $|w|\ll 1$. However, this simple procedure
which works in the case of the Yang-Mills equation mentioned above
(see \cite{BOS, RS}) does not work in the present case as we explain  below.

Note that if $\zeta \in L^2$, then $\lambda(t)$ is uniquely
determined by the orthogonality condition
\begin{equation}
\label{5a}  \langle \zeta, w\rangle =0.
\end{equation}
If $\zeta \not \in L^2$, then this condition is not well defined
unless we assume $w$ belongs to a space of sufficiently fast
decaying functions.

Substituting decomposition \eqref{5} into \eqref{2} leads to the
equation for the function $w$ and parameter $\lambda$:
\begin{equation}
\label{6}  L_xw + F(w)= - \lambda^2 \frac{\partial^2 U}{\partial
t^2},
\end{equation}
where $F(w)$ absorbs higher order terms ($F(w)= \lambda^2
\frac{\partial^2 w}{\partial t^2}+ N(w)$, $N(w)=$ nonlinearity in
$w$) and $L_x$ is the linearization operator for the r.h.s. of \eqref{1}
around $U(x)$ given by \eqref{4}. The operator $L_x$ is self-adjoint. The
scaling zero mode, $\zeta$, is a zero mode of this
operator: $L_x \zeta =0$. Since $\zeta$ is positive and not $L^2$ we
conclude by the Perron-Frobenius theory that
$\sigma(L_x)=[0,\infty)$ and $0$ is not an eigenvalue of $L_x$.

We compute explicitly
$$
 \lambda^2\frac{\partial^2 U}{\partial t^2}=\lambda^2
[-2\partial_t(\dot{\lambda}\lambda^{-1})\zeta +
2(\dot{\lambda}\lambda^{-1})^2 x \partial_{x} \zeta]
$$
\begin{equation}
\label{7}
=-2\ddot{\lambda}\lambda \zeta+2
\dot{\lambda}^2(\zeta +
 x \partial_{x} \zeta).
\end{equation}
We multiply Eqn \eqref{6} scalarly by $\zeta(x)$. Though $\zeta$ is not
$L^2$ one can show using a limiting procedure that $\langle \zeta,
L_x w \rangle = 0$, provided $w=o(x)$ and $\partial_x w=o(1)$ at
$\infty$. Thus we obtain
\begin{equation}
\label{8}
 \lambda^2\langle \zeta, \frac{\partial^2 U}{\partial t^2}
\rangle + \langle \zeta, F(w) \rangle= 0.
\end{equation}

Following \cite{BOS} we try to develop a perturbation theory in the small
parameter $\dot{\lambda}^2$ assuming that term
$\lambda\ddot{\lambda}$ is of the order $o(\dot{\lambda}^2)$ (and
$\dot{\lambda} < 0$) 
and similarly for higher
order time derivatives of $\lambda$, e.g.
$\partial_t(\lambda\ddot{\lambda})=O(\dot{\lambda}^3)$, etc. Furthermore, if our assumption that $|w|\ll 1$ is
correct and the integral in $\langle \zeta, F(w) \rangle$ is
convergent, then we can drop the term $\langle \zeta, F(w) \rangle$
in \eqref{8}. Hence we obtain in the leading order $O(\dot{\lambda}^2)$
\begin{equation}
\label{9}
 \dot{\lambda}^2\langle \zeta,\zeta +
 x \partial_{x} \zeta \rangle = 0.
 \end{equation}
Considering the integral on the l.h.s. over a bounded domain and
integrating by parts one shows that the inner product on the l.h.s.
is
\begin{equation}
\label{10}
1/2\lim_{x\rightarrow \infty}(x^2\zeta^2(x)).
\end{equation}
For Eqn \eqref{3}
this is $0$ so we can solve Eqn \eqref{6} in the leading order, $w=
-\dot{\lambda}^2 L^{-1}(\zeta +
 x \partial_{x} \zeta) $.
Plugging this result into Eqn \eqref{8} and keeping only the terms up to the
order $O(\dot{\lambda}^4)$, we obtain the equation for scaling
dynamics,
\begin{equation}
\label{11}
 \lambda\ddot \lambda=\frac{3}{4}\dot \lambda^4,
\end{equation}
in the leading order $O(\dot{\lambda}^4)$ (see \cite{BOS, RS}). Next, in order
to obtain a correction to this equation, we use \eqref{8} at the order
$O(\dot{\lambda}^6)$ to solve Eqn \eqref{6} to the order
$O(\dot{\lambda}^4)$ and plug the result to \eqref{8}. However, at this
step we run into logarithmically divergent terms. To overcome this
problem we use a multiscale expansion, by introducing an additional
scale at infinity (see \cite{BOS}).

For Eqn \eqref{1} with $k=1$ we have $\lim_{x\rightarrow \infty}(x^2\zeta^2(x))=1$
and so we go to the next term, $-2\ddot{\lambda}\lambda\|\zeta\|^2$,
and discover that it diverges logarithmically. Thus for Eqn \eqref{1} with $k=1$ one
runs into a problem right away. This shows that decomposition \eqref{5} is
incompatible with the condition $|w|\ll 1$.

The problem for Eqn \eqref{1} with $k=1$ mentioned above can be also seen in a
different but related way. 
Let us try to solve Eqn \eqref{6} by perturbation theory. In the
leading order we drop the term $F(w)$ to obtain the leading order
approximation to the solution: $w=L^{-1}\varphi$, where $\varphi  = - \lambda^2
\frac{\partial^2 U}{\partial t^2}$ and $L^{-1}$ is understood as the Green function of the equation $L w=\varphi$ (see Section \ref{approximatesolution}).  It is easy
to check, using  Eqns \eqref{12} - \eqref{13} of Section \ref{blowupvar} below, that  if $\ddot\lambda\ne
0$ then the function  $L^{-1}\varphi$ grows at infinity as $x\ln x$,
and a straightforward perturbation theory fails. (Not only the
correction $w$ is large at $\infty$, its energy is infinite.)

The point here is that the function $U(\rho/\lambda)$ is not a good
adiabatic solution to Eqn  \eqref{1} with $k=1$:
$$
\lambda^{2}(\partial_t^2U(x) - \Delta_\rho U(x) + \frac{1}{2\rho^2}\sin (2U(x)))
 $$
 \begin{equation}
 \label{14}
=-2\ddot\lambda \lambda \zeta(x)
 +\frac{4\dot\lambda^2x}{(1+x^2)^2},
\end{equation}
where $x:=\rho/\lambda$ and where we used \eqref{7} and the relation
$\zeta +
 x \partial_{x} \zeta =2x(1+x^2)^{-2}  $. The
r.h.s. is not $L^2$. The problematic term is $2\ddot\lambda \lambda \zeta(x)$. In particular, it leads to the logarithmically
divergent term, $2\ddot{\lambda}\lambda\|\zeta\|^2$  in the
orthogonality condition. Hence, one has to find a better leading
term.

We deal with the problem above by introducing instead of the linear,
one-parameter transformation, $\rho \rightarrow \rho/\lambda$, a
nonlinear, three-parameter transformation, $\rho \rightarrow
f(\rho,\lambda,\alpha, \beta )$,
chosen so that
$U(f(\rho, \lambda, \alpha, \beta))$ becomes a better
approximate solution to Eqn \eqref{1} with $k=1$ than $U(\rho/\lambda)$. In
particular, the problematic term
$2\ddot{\lambda}\lambda\zeta$ entering the r.h.s. of  Eqn
\eqref{14} 
is canceled and therefore the large  $\rho$ divergence in Eqn \eqref{8} mentioned above is eliminated.
Thus, instead of \eqref{5}, we look for solutions of Eqn \eqref{1} 
in the form
\begin{equation}
\label{15} u(\rho, t) = U(y) + w(y, t).
\end{equation}
We consider initial conditions close to $U(y) \equiv U(f(\rho, \lambda, \alpha, \beta))$ (we do not specify the norm, the latter must be determined by a rigorous analysis, see e.g. \cite{RS}).
After this we proceed as above with Eqn  \eqref{5}. 
The conditions $|w|\ll 1$ and $w\to 0$ at $\rho \to\infty$ and constraints on the energy \eqref{2_1} and its fluctuations imply the differential equation \eqref{16} on the parameter $\lambda=\lambda (t)$.
We expect that proceeding as in \cite{BOS} one can obtain corrections to Eqn \eqref{16}.

The paper is organized as follows. In Section 2 we introduce a change of
variables, $\rho \rightarrow f(\rho,\lambda,\alpha, \beta)$,  depending nonlinearly
on the original variable $\rho$ and on the scaling parameter $\lambda^{-1}$
(and depending on additional parameters $\alpha, \beta $). This is our main new
idea. In Section 3 we derive, modulo some technical details which are provided
in  Appendices 2 and 3, an approximate solution to Eqn \eqref{1} with $k=1$. In Sections 4 and 5 we use an orthogonality condition of the type of \eqref{5a}, the smallness condition on energy fluctuations and the minimum condition on the energy of the approximate solution in order to find our main equation on the scaling parameter
 $\lambda$, Eqn \eqref{16}. In Section 6 we find exact and approximate  solutions of Eqn \eqref{16} and in  Section 7 we show that this equation is  Hamiltonian. In Appendices 1-5 we provide technical calculations used in the main text and explanations of the numerical approaches.
 \DETAILS{In Appendix 1 we find two independent solutions of the equation
 $Lw=0$, see Eqns  \eqref{12} and \eqref{13}, and in Appendix 2 we
 compute some integrals occurring in the main text. In Appendix 3 we derive an
 expression for the function $\psi $, entering the transformed  equation, and in Appendix 4 we give details of computations involving the energy of the leading approximation $U(y)$.}

\section*{Acknowledgement} The authors are grateful to Piotr Bizon
for raising the wave map problem, for numerous stimulating discussions and correspondence and for providing us with numerical simulations comparing our blowup asymptotics with the ones coming from the equation. The authors are grateful to Igor Rodnianski for informing them that \cite{RRS} has succeeded to give rigorous analysis of the blowup for the problem considered here and arrived to similar conclusions as those of our non-rigorous analysis.

\section{ Nonlinear blow-up variables 
} \label{blowupvar}

In this section we introduce a nonlinear, three-parameter (scaling)
transformation, $\rho \rightarrow f(\rho,\lambda,\alpha, \beta)$ of the
independent spatial variable $\rho$. This replaces the standard,
linear, one-parameter transformation, $\rho \rightarrow
\rho/\lambda$. We denote  $x:=\rho/\lambda$ and
define $y= y(x,\lambda\ddot\lambda,\alpha, \beta)$ as the solution of the equation
$$
 y = x - \frac{\lambda\ddot\lambda}{2}x^{3}
\ln\biggl(\frac{\sqrt{\lambda\ddot\lambda}}{ \beta}
y^{\alpha} x^{1-\alpha}\biggr),\qquad\mbox{ if}\   x\le x_{cr} $$
\begin{equation}
\label{20} y = 2y_{cr} -x + \frac{\lambda\ddot\lambda}{2} x^{3} \ln\biggl(
\frac{\sqrt{\lambda\ddot\lambda}}{ \beta}y^{\alpha}x^{1-\alpha}\biggr),\qquad\mbox{if }\ x>x_{cr},
\end{equation}
where $0 \le \alpha \le 1,\ \beta >0$, $x_{cr}=x_{cr}(\lambda, \alpha, \beta)$ 
and $y_{cr}= y_{cr}(\lambda,\alpha, \beta)$
solve the equations $\frac{\partial }{\partial x}(\mbox{r.h.s. \eqref{20}}) = 0$ and the first equation in \eqref{20}.
We can write $x_{cr}$ and $y_{cr}$ 
as
\begin{equation}
\label{22}   x_{cr}= \biggl(\frac{2}{\lambda\ddot\lambda}\biggr)^{1/2}  (
3~\ln\gamma +1-\alpha)^{-1/2}
\end{equation}
and
\begin{equation}
\label{23}   y_{cr} = \biggl(\frac{2}{\lambda\ddot\lambda}\biggr)^{1/2}
\frac{2\ln\gamma+1-\alpha}{(3\ln\gamma+1-\alpha)^{3/2}},
\end{equation}
where $\gamma = \gamma (\alpha, \beta)=\frac{\sqrt{\lambda\ddot\lambda}}{ \beta}y^{\alpha}x^{1-\alpha}$ is a solution of the equation
\begin{equation}
\label{24}  \frac{ \gamma}{\sqrt{2}}\frac{(3\ln\gamma+1-\alpha)^{1/2+\alpha}} {(2\ln\gamma+1-\alpha)^{\alpha}}=\beta^{-1}.
\end{equation}
 Eqn \eqref{24} is well defined for $\ln\gamma > -\frac{ 1}{3}(1-\alpha)$ and in this range it has a unique solution for each $\alpha$ and $\beta$. We denote  $\rho_{cr}:= \lambda x_{cr}$.

For $x \le x_{cr}$ the r.h.s. of \eqref{20} is decreasing from $\infty$ to $-\infty$. 
Hence the equation \eqref{20} has a unique solution for $x \le x_{cr}$. Furthermore, for $x \ge x_{cr}$  the r.h.s. of \eqref{20} increases logarithmically  in $y$ from $-\infty$ to $\infty$ and 
 increases monotonically in $x$. Since for $x=x_{cr}$ \eqref{20} has a unique solution, it has exactly two solutions  for $x > x_{cr}$. Of these two solutions we choose the greater one.

 Finally, we notice that the  function $y= y(x,\lambda\ddot\lambda,\alpha, \beta)$ increases monotonically in $x$ for $x >0 $. Indeed,   for $  0 \le\alpha \le 1 $, the equations $\frac{\partial }{\partial x}(\mbox{r.h.s. \eqref{20}}) =0$ and \eqref{20} have a unique solution ($x=x_{cr},\ y=y_{cr}$), and therefore so is the equation $\frac{\partial y}{\partial x} =0$. Since $y$ is obviously  increases monotonically in $x$ for $x  $ sufficiently small and sufficiently large, it does so for all  $x >0 $.

Write $v(y,t)=u(\rho,t)$, where $y=y(x,\lambda\ddot\lambda,\alpha, \beta)$ is given in
Eqn \eqref{20}. In the new variables, Eqn \eqref{1} with $k=1$ becomes
\begin{equation}
\label{25}
-\frac{\partial^2v}{\partial y^2} -
\frac{1}{y}~\frac{\partial v}{\partial y} + \frac {\sin (2v)}{2y^2}
=\Psi(v),
\end{equation}
where
\begin{equation}
\label{26} \Psi(v):= \frac{x^2}{y^2} \biggl \{ \left(\frac{2y}{x }\chi
+\chi^2\right)\frac{\partial^2v}{\partial y^2}+ \left(\frac{2\chi}{x}
+\frac{\partial\chi}{\partial x}\right)\frac{\partial v}{\partial y} - \lambda^2 \frac{\partial^2v}{\partial
t^2}\biggr\}
\end{equation}
 and $\partial^2/\partial t^2$ is the total derivative in
$t$ (i.e. taking into account that $y$ depends on $t$).
Here the function $\chi$ is defined according to the equation
\begin{equation}
\label{27'} \frac{\partial y}{\partial x} = \frac{y}{x} + \chi.
\end{equation}
Eqn \eqref{25} is our transformed equation.

Initial conditions for \eqref{25} are chosen to be close, in an appropriate norm, to $U(y) \equiv U(y(x, \lambda\ddot\lambda, \alpha, \beta))$, where, recall, $U(\rho)$ is the static - kink - solution
to  Eqn \eqref{1}. To simplify the exposition we take the initial condition to be just  $U(y(x, \lambda\ddot\lambda, \alpha, \beta))$.

\section{ Approximate solution of Eqn \eqref{25}} \label{approximatesolution}

Let $y= y(x,\lambda\ddot\lambda,\alpha, \beta)$ be the transformation defined in the previous section. We look for solutions of Eqn \eqref{25} in the form
\begin{equation}
\label{26a} v(y, t) = U(y) + w(y, t),
\end{equation}
where $w$ is a small correction. We plug this decomposition into Eqn
\eqref{25} to obtain \begin{equation}
\label{27}
L_y w +N(w)= \Psi(U+w),
\end{equation} where operator $L_y$ is
defined in Eqn \eqref{4}, $N(w)$ is the nonlinear in $w$ term defined by this equation and the function $\Psi(v)$ is defined in \eqref{26}. To find an approximate solution of the latter equation we drop the nonlinearity, $N(w)$, and the term $w$ in
$\Psi(U+w)$ to obtain the leading order equation
\begin{equation}
\label{28}  L_y W =  \psi,
\end{equation}
where $\psi(y, t):=\psi(U(y)) $. The latter function is given explicitly by
\begin{multline}
\label{29}
 \psi(y, t) := \frac{x^2}{y^2} \biggl \{ \frac{8\chi}{x
(1+y^2)^2} + \frac{2}{1+y^2} \biggl ( \frac{\partial\chi}{\partial x} - \frac{2\chi}{x} \biggr ) -
\frac{4y\chi^2}{(1+y^2)^2} \\ + \frac{4 \lambda^2y}{(1+y^2)^2}\biggl(\frac{\partial y}{\partial t}\biggr)^2
-\frac{2 \lambda^2}{1+y^2}\frac{\partial^2 y}{\partial t^2}\biggr \},
\end{multline}
where the variable $x$ is connected to $y$ through \eqref{20}. The counter-term which removes the undesirable term in $\frac{\partial^2 U}{\partial t^2}$ is
$\partial\chi/\partial x - 2\chi/x$. As the result we have
\begin{equation} \label{asbehavpsi}\psi(y, t) =O\biggl(\dot\lambda^2\frac{\ln (\lambda\ddot\lambda y^2)}{\lambda\ddot\lambda
y^2}\biggr)\ \mbox{for}\ y \ge y_{cr}
\end{equation}
(see Appendix 3).
\begin{remark} To justify dropping $\frac{\partial^2 w}{\partial t^2}$ from Eqn \eqref{27} (it enters into $\psi(U+w)$) is not a simple matter. In a rigorous approach one looks for solutions of \eqref{27} in the form
\begin{equation*}
 v(y, t) = U(y) + W(y)+\eta(y, t)
\end{equation*}
and shows that the fluctuation $\eta$ is small (cf. \cite{RS}).
\end{remark}

Recall that the  operator $L_y$ entering  Eqn \eqref{28} has a zero mode:
\begin{equation}
\label{35}
 L_y \zeta  = 0.
\end{equation}
Multiplying Eqn \eqref{28} scalarly by $\zeta$ and using  the
self-adjointness of the operator $L$ (and some elementary limiting
procedure) and \eqref{35}, we obtain
\begin{equation}
\label{36} \int\limits^{\infty}_0 dy y\ \zeta \psi = 0.
\end{equation}
This is a (necessary) solvability condition for  Eqn \eqref{28}. It gives an equation on the parameters  $\lambda,\ \alpha$ and $\beta$. (Note that it is an approximate solvability condition for the exact equation \eqref{27}.

So far we obtained one equation for the three parameters  $\lambda,\ \alpha$ and $\beta$. To derive another equation we analyze the approximate solution to \eqref{27} we obtained: $U(y)+W(y)$, where $W=L^{-1}\psi$ 
with $\psi$ satisfying \eqref{36}. 
Our goal in the rest of this section is to isolate the leading contribution to $W$. This will be used in the next section to derive the second equation for the parameters.

To find $L^{-1}\psi$ we compute the Green's
function for the operator $L$. Two linearly independent solutions of
the homogeneous equation $L w = 0$ are
\begin{equation}
\label{12}
 w_1(y) = \frac{y}{1+y^2}  \qquad  \mbox {\rm and} \qquad w_2(y)
= \frac y2 - \frac{1}{2y} + \frac{2y\ln y}{1+y^2}
\end{equation}
(the first of these solutions is just the scaling zero mode,
$\zeta$, the second solution is found in Appendix 1). Hence by the
ODE theory
\begin{equation}
\label{13}
 L^{-1}\psi = cw_1+ w_1\int_{0}^{y}w_2\psi y'dy'  -
w_2\int_{0}^{y}w_1\psi y'dy'
\end{equation}
where $c$ is chosen to guarantee solvability of the equation to the
second order correction term or by minimizing the energy.

We find the leading contribution to the solution $w=L^{-1}\psi$ of  Eqn \eqref{28}.
In what follows we use the following assumptions
\begin{equation}
 \label{33} 0<\lambda\ddot\lambda \ll \dot\lambda^2 \ll 1,\ \lambda \partial_t ( \lambda\ddot\lambda ) = O(\dot\lambda^3),\
\lambda \partial_t ( \dot\lambda^2) = O(\dot\lambda^3),\  \beta =O(1).
 \end{equation}
We see from \eqref{22} that $$x_{cr}=\rho_{cr}/\lambda\sim (\lambda\ddot\lambda)^{-1/2} \gg 1.$$

Consider first the region $y \le y_{cr}$. In this region $\psi$ is given by \eqref{A3.1}, Appendix 3. The latter equation shows that for $y \ll y_{cr}$ the leading part of $\psi$ is
\begin{equation}
\label{38} \psi_1 (y, t) := -\frac{x^3}{y^2} \frac{8\lambda\ddot\lambda
\left(\ln\left(\frac{\sqrt{\lambda\ddot\lambda}}{\beta}y^\alpha x^{1-\alpha}\right)+1/2\right)  - 4\dot\lambda^2}{(1+y^2)^2\left(1 +
\frac{\alpha\lambda\ddot\lambda x^3}{2 y}\right)},
\end{equation}
where, recall, $x:= \rho/ \lambda$ is connected to $y$ through \eqref{20}. Now, let
\begin{equation}
\label{39}
\psi_2:=\psi - \psi_1.
\end{equation}

According to \eqref{13}, the general solution, $W$, of Eqn \eqref{28} decreasing at infinity is of the form
\begin{equation*}
 W(y, t) = c_1 w_1(y) + w_1(y)\int_{0}^{y}w_2(s)\psi(s,t)
sds
\end{equation*}
$$ + w_2(y)\int_{y}^{\infty}w_1(s) \psi(s,t) sds,
$$
with $w_1=\zeta$ and $w_2$ defined in \eqref{12}. The function $w_2$ is
singular at $y=0$. Hence $W$ is bounded only if the condition \eqref{36}
is satisfied.

The function $\psi_1(y, t)$ is localized on the scale $y\sim 1$, decays as
$y^{-2}$ in the region  $1 \ll y \ll y_{cr}$ and decays as $y^{-5}$
in the region $ y \gg y_{cr}$ (though we are considering at the
moment only  the region $1 \ll y \ll y_{cr}$, the latter fact allows us to
extend the integrals to the entire real axis with a small error).
The function  $\psi_2$ is localized at the large scale, $ y\sim y_{cr} \gg 1$. After some lengthy computations we find for $y \ll  y_{cr}$
\begin{equation}
\label{42}
 W(y,t) = \frac{y}{2}
\int_{0}^{\infty}w_1(s)\psi_1(s,t) s ds
+O(\dot\lambda^2\lambda\ddot\lambda{y}^3).
\end{equation}
Using the expression \eqref{38}, it is easy to show that  $$\int_{0}^{\infty}w_1(s)\psi_1(s,t)sds=O(\dot\lambda^2).$$

\section{ Energy of the approximate solution and the equation on $\lambda$}

We compute the energy of our approximate solution $u(\rho, t)=U_{\lambda, \alpha, \beta}(\rho) + W(y)$, where $U_{\lambda, \alpha, \beta}(\rho)  :=
U(y)$, with $y=y(x, \lambda\ddot\lambda, \alpha, \beta)$ and $U$ defined in the
introduction, and  where $W=L^{-1}\psi$, the solution to Eqn \eqref{28} (see the previous section).  
Due to \eqref{2_1}, the energy functional is
\begin{equation}
 \label{37}
 E(u) :=\int_0^\infty (\frac{1}{2} \dot{ u}^2 + \frac{1}{2}|\nabla u|^2 +  \frac{1}{2\rho^2}\sin^2u)\rho d\rho.
  \end{equation}
Inserting the approximate solution into this expression, we obtain that $E(U_{\lambda, \alpha, \beta} + W) = E(U_{\lambda, \alpha, \beta}) + \delta E_1$ with
\begin{equation}
\label{43}
 \delta E_1 :=O(\dot\lambda^2y_{cr}^4)
 \biggl(\int_{0}^{\infty}w_1(s)\psi_1(s,t) s ds\biggl)^2+O(\dot\lambda^2).
\end{equation}
Furthermore, we have that
 \begin{equation}
\label{deltaE0}
E(U_{\lambda, \alpha, \beta}) = E(U) + \delta E_0\ \mbox{with}\  \delta E_0 = O(\dot\lambda^2 \ln(1/\dot\lambda^2)).
\end{equation}
 We require that the energy
correction due to the fluctuation, $W$, be much smaller than the one due to the modulation:
\begin{equation}
\label{deltaE0detalE1}
|\delta E_1| \ll |\delta E_0|.
\end{equation}
 Since $\int_{0}^{\infty}w_1(s)\psi_1(s,t)sds=O(\dot\lambda^2)$ and $y_{cr}= O(\frac{1}{ \sqrt{\lambda\ddot\lambda}})\gg O(\frac{1}{ |\dot\lambda|})$, this implies that the integral in the
leading term in the above expression for $\delta E_1$ must vanish:
\begin{equation}
\label{44}\int\limits^{\infty}_0  dyy \zeta \psi_1 = 0.
\end{equation}
This gives an implicit equation on the parameters $\lambda$, $\alpha, \beta$.

In the leading order, we can replace $y$ by $x=\rho/\lambda$ (see the first equation in \eqref{20}), so that Eqn \eqref{44}   becomes
\begin{equation}
\label{45} \int\limits^{\infty}_0 dx x \frac{x}{1+x^2} \biggl \{ \frac{8\lambda\ddot\lambda x}{(1+x^2)^2}
\biggl [ \ln \left(\frac{\sqrt{\lambda\ddot\lambda}}{ \beta}x\right) + 1/2 \biggr ] + \frac{4\dot\lambda^2
x}{(1+x^2)^2} \biggr \} = 0.
\end{equation}
Computing the integrals in \eqref{45} (see Appendix 2 for detailed computations), we obtain
\begin{equation}
\label{46} \dot\lambda^2 + 2\lambda\ddot\lambda \biggl [
~\ln\left(\frac{\sqrt{\lambda\ddot\lambda}}{ \beta}\right) +1 \biggr ] = 0.
\end{equation}
This is our explicit
equation for the parameter $\lambda$. It depends on the additional parameter $ \beta$ whose value we still have to determine. Since in the leading approximation ($y \rightarrow x$) the first equation on the r.h.s. of \eqref{20} is independent of $\alpha$, then so are the resulting equations \eqref{45} and \eqref{46}.
 Eqns  \eqref{46} and  \eqref{16} coincide, provided
\begin{equation}
\label{49}a =  \beta^2\ e^{-2}.
\end{equation}

Clearly, solutions of Eqn \eqref{46} 
have the property \eqref{33} assumed above.
Moreover, if $\lambda(0) >0,\ \dot\lambda(0)<0 
$, then,  by Eqn \eqref{46}, 
$\lambda(t) >0,\ \dot\lambda(t)<0,\ \ddot\lambda (t)>0$ and therefore $ \dot\lambda(t)^2 \le \dot\lambda(0)^2 
$ for $t >0$.
%
As $t\to t_*,\ |\dot\lambda|$ decreases so that our approximation improves as $t\to
t_*$.
%

Thus it remains to find the value of the parameter $\beta$ .
To this end we use the condition \eqref{36} and minimization of the energy of the leading part of the approximate solution.

\section{Values of the parameters $\alpha$ and $\beta$ }
\label{sec:4}

\bigskip

In this section we derive an equation on the parameters  $\alpha$ and $\beta$ and use this equation together with the energy minimization to find the values of these parameters.
We assume \eqref{33} and that, at least in the leading approximation,
\begin{equation}
\label{46b}\alpha, \beta\ \mbox{are independent of}\ t.
\end{equation}

\textit{In what follows we do not display the dependence of the quantities involved on $\lambda\ddot\lambda$ (however, the dependence on $\lambda$ is displayed)}.

First note that Eqns \eqref{36},  \eqref{39} and \eqref{44} imply that
\begin{equation}
\label{52}
\int_{0}^{y_{cr}}\psi_{2}dy+\int_{y_{cr}}^{\infty} \psi dy=0
\end{equation}
From Eqns \eqref{23} and \eqref{46} we obtain easily
$$
\frac{1}{\lambda\ddot\lambda}\lambda\frac{\partial(\lambda\ddot\lambda)}{\partial t} =
2 \dot\lambda \ln^{-1}\left(\frac{1}{\lambda\ddot\lambda}\right)\ \left[1+O\left(1/\ln\left(\frac{1}{\lambda\ddot\lambda}\right)\right)\right]
$$
and
\begin{equation}
\label{53}
\frac{1}{x_{cr}}\lambda\frac{\partial x_{cr}}{\partial t} = -\dot\lambda\ln^{-1}\left(\frac{1}{\lambda\ddot\lambda}\right)\ \left[1 +
O\left(1/\ln\left(\frac{1}{\lambda\ddot\lambda}\right)\right)\right].
\end{equation}

As a result in the main approximation in $1/\ln(1/\lambda\ddot\lambda)$ we should
keep in the expressions for the functions $\psi$ and $ \psi_2$ in Eqn \eqref{52} only terms proportional to $\dot\lambda^2$. The latter terms are given in  \eqref{A3.1} and \eqref{A3.2} in Appendix 3. The most important region in the above integral is where $y$ is of order of $y_{cr}$ ($y_{cr}=O(\frac{1}{ \sqrt{\lambda\ddot\lambda}})\gg 1$ due to Eqn \eqref{23} and the condition \eqref{33}). As a result we can neglect
$1$ compared to $y^2$ in \eqref{A3.1} and \eqref{A3.2}.  Then Eqn \eqref{52}  reduces to the equation
\begin{equation}
\label{54}
I( \alpha, \beta): = I_1( \alpha, \beta) + I_2( \alpha, \beta) =0,
\end{equation}
where
\begin{multline}
\label{55}
I_1 =2 \int\limits_0^{y_{cr}}\frac{dy x^4}{y^4}\biggl\{\frac{\lambda\ddot\lambda x^2
\left(\ln\left(\frac{\sqrt{\lambda\ddot\lambda}}{\beta}y^\alpha x^{1-\alpha}\right) -
\alpha\right)}{y\left(1 + \frac{\alpha\lambda\ddot\lambda x^3}{2 y}\right)^2}   -
\frac{2(1-Y^2)}{y\left(1 + \frac{\alpha\lambda\ddot\lambda x^3}{2 y}\right)^2} +\frac{2(1-Y)}{x(1 +\frac{\alpha\lambda\ddot\lambda x^3}{2 y})}  \\-
\frac{\lambda\ddot\lambda x}{2} \biggl(\frac{\alpha x^2}{y^2}\frac{Y^2}{\left(1 +
\frac{\alpha\lambda\ddot\lambda x^3}{2 y}\right)^3} - \frac{6\alpha x}{ y\left(1 +
\frac{\alpha\lambda\ddot\lambda x^3}{2 y}\right)^2} - \frac{3(1-\alpha)}{1 +
\frac{\alpha\lambda\ddot\lambda x^3}{2 y}} \\
- \frac{2\left(1 -\frac{\alpha\lambda\ddot\lambda x^3}{y}\right)\left(3\ln\left(\frac{\sqrt{\lambda\ddot\lambda}}{\beta}y^\alpha x^{1-\alpha}\right) + 1 -
\alpha\right)}{\left(1 + \frac{\alpha\lambda\ddot\lambda x^3}{2 y}\right)^2}\biggr)\biggr\}
\end{multline}
and
\begin{equation}
\label{55_1}
I_2 =2 \int\limits_{y_{cr}}^\infty dy\frac{x^4}{y^4}\biggl\{\frac{2Y^2}{y\left(1
- \frac{\alpha\lambda\ddot\lambda x^3}{2 y}\right)^2} +
\frac{2Y}{x(1 - \frac{\alpha\lambda\ddot\lambda x^3}{2 y})}
 $$$$ + \frac{\lambda\ddot\lambda x}{2}
\biggl[\frac{\alpha x^2}{y^2}\frac{Y^2}{\left(1 -
\frac{\alpha\lambda\ddot\lambda x^3}{2 y}\right)^3} + \frac{6\alpha x}{ y\left(1 -
\frac{\alpha\lambda\ddot\lambda x^3}{2 y}\right)^2} - \frac{3(1-\alpha)}{1 -
\frac{\alpha\lambda\ddot\lambda x^3}{2 y}} $$$$  - \frac{2\left(1 +
\frac{\alpha\lambda\ddot\lambda x^3}{
y}\right)\left(3\ln\left(\frac{\sqrt{\lambda\ddot\lambda}}{\beta}y^\alpha x^{1-\alpha}\right) + 1 -
\alpha\right)}{\left(1 - \frac{\alpha\lambda\ddot\lambda x^3}{2 y}\right)^2}
\biggr]\biggr\}.
\end{equation}
Here
\begin{equation}
\label{30_2'}
Y=1-\frac{\lambda\ddot\lambda x^2}{2}\left(3\ln\left(\frac{\sqrt{\lambda\ddot\lambda}}{ \beta}
y^\alpha x^{1-\alpha}\right)+1-\alpha\right).
\end{equation}

One can further evaluate $I_1$ and $I_2$ by changing the variable of integration $y$ in \eqref{55} and \eqref{55_1} to $z$ as
\begin{equation}
\label{56}
z=x/x_{cr} =\rho / \rho_{cr}, 
\end{equation}
where $x$ and $x_{cr}$, as functions of $y,\ \lambda\ddot\lambda,\ \alpha,\ \beta$ are given in \eqref{20} and the definitions following this equation,
and compute the resulting integral numerically.   In particular, one can show
that for $\alpha = 0, I_1( \alpha=0, \beta) = 1, I_2( \alpha=0, \beta) = 0$ and therefore  $ I( \alpha=0, \beta) = 1$, independently of the value of $\beta$. Thus we cannot take $\alpha=0$ in our transformation \eqref{20}.

We chose the parameters $\alpha$ and $\beta$ which minimize the energy $E( \alpha, \beta):= E(U_{\lambda, \alpha , \beta})$, where, recall,  $U_{\lambda, \alpha, \beta}(\rho) :=U(y(x, \lambda\ddot\lambda, \alpha, \beta))$, given that the equation  \eqref{54}, $I(\alpha, \beta)=0$, holds.
To find these minimizers we use Eqns \eqref{37} and $U(\rho) = 2~\mbox {\rm  arctan}~ \rho$ to rewrite the energy $E( \alpha, \beta)$ as
\begin{equation}
\label{57}
E (\alpha, \beta)= 2\int\limits_0^\infty d\rho\rho\frac{1}{(1+y^2)^2}\left\{\left(\frac{\partial y}{\partial t}\right)^2 +
\left(\frac{\partial y}{\partial\rho}\right)^2 + \frac{y^2}{\rho^2}\right\}.
\end{equation}
We find numerically (see Appendix 5 for the analytical part) that the energy $E( \alpha, \beta)$ is minimized on the curve $I(\alpha, \beta)=0$ at the point
\begin{equation}
\label{61_2}
\beta_0 = 1.04\ \mbox{and}\ \alpha_0 = 0.65436.
\end{equation}

This is a special point for the curve $I(\alpha, \beta)=0$. Our numerics show that while the functions  $\alpha = \alpha( \beta)$ and  $\beta = \beta( \alpha)$ determined by the equation $I(\alpha, \beta)=0$ are double-valued, their branches originate exactly at this point (and form a wedge there). So the equation $I(\alpha, \beta)=0$ has a unique solution only for $\beta=\beta_0$ or for $\alpha=\alpha_0$ and has no solutions for  $\beta> \beta_0$ or for $\alpha<\alpha_0$.

Substituting $\beta = \beta_0= 1.04$ into Eqn \eqref{46},
we obtain the following value for the parameter $a$:
$$ 
a = 0.146.$$
This proves Eqn  \eqref{16} with $a = 0.146.$


\bigskip


\section{ Investigation of Eqn \eqref{16}} \label{eqn2}

In this section we find an approximate solution to Eqn \eqref{16} (which is, up to a redefinition of the parameters, the equation \eqref{46}).
Iterating this equation, we find, in the leading approximation, the
following equation
\begin{equation}
\label{65}
 \frac{\ddot\lambda}{\dot\lambda} =
\frac{\dot\lambda}{\lambda \ln ( a/\dot\lambda^2 )}.
\end{equation}
Solution of the Eqn \eqref{65} with two free parameters of
integration, $ c>0$ and $t^{*}$, is
\begin{equation}
\label{66}
 \sqrt{a}~(t^{*} - t) = \int\limits^{\lambda}_0~dx~e^{ \ln^{1/2}  ( \frac {c}{x}  ) }.
\end{equation}
Changing the variable of integration as
$\ln \biggl( \frac{c}{x}\biggr) = z^2,$
we reduce Eqn \eqref{66} for the parameter $\lambda$ to the form Eqn \eqref{17} given in Introduction.


Now we derive an exact expression for a general solution of \eqref{16}. We introduce the function $f(x) : = x\ln (1/x)$. For $0 < x <
e^{-1}$ this function has the inverse, $f^{-1} (x)$. Using this
inverse we rewrite Eqn \eqref{16} as

\begin{equation}
\label{69}
 \frac{\lambda \ddot\lambda}{a} = f^{-1} \biggl({\dot{\lambda}^2
\over a}\biggr).
\end{equation}
(Note that for $x \rightarrow 0,\ f^{-1}(x)=\frac{x}{\ln x}+...$, so in the leading  approximation of \eqref{69} gives \eqref{65}.) Integrating equation \eqref{69} gives
\begin{equation}
\label{70}
 \ln \lambda = F (\dot{\lambda})\ \textrm{where}\ F(y) = {1
\over 2} \int\limits^{{y^2 \over a}} {dz \over z} g(z),
\end{equation}
 with $g(z) : = z /{f^{-1}(z)}$.  Using the equation
$f(f^{-1} (y)) = y$, or, more explicitly, $f^{-1}(y)\ln (1/ {f^{-1}(y)})
= y$, we find that the function $g(z)$ satisfies the equation
\begin{equation}
\label{71}
 g(z) = \ln\biggl({g(z) \over z}\biggr).
\end{equation}
\medskip
 Differentiating the latter equation, we find
 \begin{equation}
 \label{72}
 g'(z) = -\frac{g(z)}{z(g(z)-1)}.
 \end{equation}
 Using this equation we integrate
\begin{equation}
\label{73}
 \int\limits^x {dz \over z} g(z) = - \int\limits^x dz
g'(z)(g(z) -1) = - {1 \over 2} (g(x) -1)^2 + \textrm{const}.
\end{equation}
 This gives
 \begin{equation}
 \label{74}
 F(y) = - {1 \over 4}(g ({y^2/ a})
-1)^2 + const,
\end{equation}
 which together with Eqn \eqref{70} yields
\begin{equation}
\label{75}
 g\biggl({\dot{\lambda}^2 \over a}\biggr) = 1 + 2 \sqrt{\ln\biggl({c \over
\lambda}\biggr)}
\end{equation}
for some constant $c$. The latter equation can be integrated as
follows
\begin{equation}
\label{76}
 \sqrt{a} (t^* - t) = \int_0^\lambda\frac{dx}{
\biggl[g^{-1}\biggl(1 + 2 \sqrt{\ln (\frac{c}{x}) }\biggr)\biggr]^{1/2}}.
\end{equation}
 Next we find the function
$g^{-1} (x)$. The definition of the function $f(x)$ implies
$f(e^{-x}) = x e^{-x}$, which yields
\begin{equation}
\label{77}
 \frac{x e^{-x}}{f^{-1}(x e^{-x})} = x ,
\end{equation}
which, in turn, leads to $g(x e^{-x}) = x$, which finally gives the
expression
\begin{equation}
\label{78}
 g^{-1} (x) = x e^{-x}.
\end{equation}
 Now Eqns \eqref{76} and \eqref{78} imply
\begin{equation}
\label{79}
 \sqrt{a} (t^* - t) = \int\limits_0^{\lambda} dx \frac{e^{1/2+
\sqrt{\ln(\frac{c}{x})}}}{\sqrt{1 + 2
\sqrt{\ln(\frac{c}{ x})}}}.
\end{equation}
 Changing the variable of integration as $\ln (c/x)  = z^2$
we find
\begin{equation}
\label{80}
 \sqrt{a}(t^*-t) = \sqrt{2}ce^{1/2}\int\limits^\infty_{\sqrt{\ln(c/\lambda)}} dz
 \frac{z}{\sqrt{z+1/2}}
 e^{z-z^2}.
\end{equation}
We obtain  for $\sqrt{\ln (c /\lambda)} \gg 1$ the
approximate expression
\begin{multline}
\label{81}
 \sqrt{a}(t^*-t)=\frac{\lambda}{\sqrt{2}[\ln(c/\lambda)]^{1/4}}
 e^{1/2+\sqrt{\ln(c/\lambda)}}\\\times
 \biggl[1+\frac{1}{4\sqrt{\ln(c/\lambda)}}-\frac{1}{32}\frac{1}{\ln(c/\lambda)}\biggr].
\end{multline}
Eqn \eqref{17} is an approximation for this exact expression, it differs from the latter by a slowly varying factor which can be found in the next approximation to \eqref{17}. 

\section{ Hamiltonian Formulation} \label{Ham}

\medskip
 Eqn \eqref{16} is a Hamiltonian system.  Indeed, it can be obtained from
the Langrangian
\begin{equation}
\label{82}
  L= h (\dot\lambda) - \ln \lambda
\end{equation}
 where the function $h$ is defined by
\begin{equation}
\label{83}
  h'' (x) = -\frac{1}{f^{-1}(x^2/a)}
\end{equation}
 with $f (x) = x \ln ( 1/ x)$ (see Section 5). Now the generalized momentum,
Hamiltonian and energy can be found in the standard way. In
particular, the energy is given by
\begin{equation}
\label{84} E =-\dot \lambda {\partial L \over \partial \dot\lambda}
+L = -\dot\lambda h' (\dot\lambda) + h (\dot\lambda)-\ln \lambda.
\end{equation}
This is the energy
conservation law.  On the other hand, differentiating Eqn \eqref{84} w.r.to $t$,  we obtain the equation of motion  \eqref{16}.

\section*{7. Conclusion}

\bigskip

We presented detailed arguments that for an open set of initial conditions close to the degree $1$ equivariant, static wave map, the solutions of the wave map equation ($\sigma$-model) collapse in a finite time. Near the collapse point the solutions have a universal profile given by the modified  the degree $1$ equivariant, static wave map depending on a  time-dependent parameter $\lambda$.
This parameter describes the rate of compression (scaling) of the collapse profile. We derived a second order Hamiltonian dynamical equation for the scaling parameter, $\lambda$. We also found approximate solutions of this equation.  These solutions are of a rather complex form. They are in an excellent agreement with direct numerical simulations of the wave map equation. 

\bigskip




\bigskip

\textbf{Appendix} 1

\bigskip

To solve the equation $Lw=g$, we should find first of all two linear
independent solution of linear  equation
\begin{equation}
 L_xw = 0. \label{A1.1}
 \end{equation}
The first solution of this equation is the scaling zero mode
$\zeta$
\begin{equation}
w_1= \zeta=\frac{x}{1+x^2}.
\label{A1.1}
 \end{equation}
The second solution $w_2$ satisfies the inhomogeneous  equation of
first order:
\begin{equation}
w_1w_2' - w_2w_1' = \frac {1}{x}. \label{A1.2}
 \end{equation} The standard
solution of this equation is
\begin{equation}
w_2 = w_1z; \quad z' = x + \frac {2}{x} +\frac{1}{x^3}; \quad z = C
+ \frac{x^2}{2} + 2\ln x - \frac {1}{2x^2}.  \label{A1.3}
 \end{equation}
Setting
$C=0$, we obtain
\begin{equation}
w_2 =\frac {x}{2} + \frac{2x\ln x}{1+x^2}- \frac {1}{2x}. \label{A1.4}
 \end{equation}

To obtain general solution of the equation $L_xw=g$, we
rewrite it as a first order ODE
\begin{equation}
\frac{\partial}{\partial x}  {w\choose v} = \biggl (
\begin{array}{ll}
0  & 1 \\
\frac{1}{x^2}\biggl ( 1-\frac{8x^2}{(1+x^2)^2} \biggr ) & - \frac
{1}{x}
\end{array}
\Biggr ) {w\choose v} - g {0\choose 1}.
 \label{A1.5}
 \end{equation}
Two linear independent solutions of A1.5 are
 \begin{equation}
 {w_1\choose w_1'},  \quad  {w_2\choose w_2'}.  \label{A1.6}
 \end{equation}
By the method of variation of constants we look for a general
solution of inhomogeneous Eqn \eqref{A1.5} in the form
  \begin{equation}
  {w\choose v} = c_1 {w_1\choose w_1'} + c_2 {w_2\choose w_2'}.
  \label{A1.7}
 \end{equation}
where $c_{1,2}$  are functions of $x$. Inserting  \eqref{A1.7} into Eqn
\eqref{A1.5}, we find
\begin{equation}
\frac{\partial c_1}{\partial x} = xw_2g;  \qquad
\frac{\partial c_2}{\partial x} =-xw_1 g. \label{A1.8}
 \end{equation}

\bigskip

\textbf{Appendix} 2

\bigskip

To derive Eqn \eqref{46} from Eqn \eqref{45} we should calculate two simple
 integrals. One of them is
 \begin{equation}
 \int\limits^{\infty}_0\frac{dy~y^3}{(1+y^2)^3} = \frac {1}{2}
 \int\limits^{\infty}_0 ~\frac{dx~x}{(1+x)^3} = \frac {1}{4}.
  \label{A2.1}
 \end{equation}
 The second integral  is $(\varepsilon\to 0)$
 \begin{equation}
\int\limits^{\infty}_0 \frac{dy~y^3\ln y}{(1+y^2)^3} = \frac {1}{4}
\int\limits^{\infty}_0 \frac{dx~x\ln x}{(1+x)^3} = -\frac {1}{4}
\int\limits^{\infty}_{\varepsilon} \ln x~d \biggl ( \frac{1}{x+1} -
\frac {1}{2}~ \frac{1}{(x+1)^2}  \biggr )     \label{A2.2}
 \end{equation}
\begin{equation*}
= \frac {1}{8} \ln\varepsilon+\frac {1}{4}
\int\limits^{\infty}_{\varepsilon}~\frac {dx}{x} \biggl (
\frac{1}{x+1} - \frac {1}{2}~ \frac{1}{(x+1)^2} \biggr )
\end{equation*}
\begin{equation*}
=\frac {1}{8} \ln\varepsilon + \frac {1}{8}
\int\limits^{\infty}_{\varepsilon}~dx \biggl [ \frac {1}{x} - \frac
{1}{x+1}+ \frac {1}{(x+1)^2} \biggr ] = \frac {1}{8}.
\end{equation*}
Using the values of these two integrals, we obtain Eqn \eqref{46}  from Eqn
\eqref{45}.

\bigskip

\textbf{Appendix} 3

\bigskip

Now we will find an explicit expression for the inhomogeneous term $\psi$. We consider separately two domains $\{y
\le y_{cr}\} \equiv \{\rho \le \rho_{cr}\}$ and $\{y \ge y_{cr}\} \equiv \{\rho \ge \rho_{cr}\}$. First, we compute  $\frac{\partial y}{\partial t}$ and $\frac{\partial^2y}{\partial t^2}$.

 Recall the notation  $x:=\rho/\lambda$ and
\begin{equation}
\label{30_2}
Y:=1-\frac{\lambda\ddot\lambda x^2}{2}\left(3\ln\left(\frac{\sqrt{\lambda\ddot\lambda}}{ \beta}
y^\alpha x^{1-\alpha}\right)+1-\alpha\right)
\end{equation}
and let
\begin{equation}
\label{XZ}
A:= \ln\left(\frac{\sqrt{\lambda\ddot\lambda}}{ \beta}
y^\alpha x^{1-\alpha}\right),\ X:=1 +
\frac{\alpha\lambda\ddot\lambda x^3}{2y} \ \mbox{and}\ Z:=1 - \frac{\alpha\lambda\ddot\lambda x^3}{2 y}
\end{equation}
In the domain $y \le y_{cr}$ we have
$$
\lambda\frac{\partial y}{\partial t}=-\dot\lambda x\frac{Y}{X} +
O
    \left(
        x^3\lambda\frac{\partial}{\partial t}\left(
            \lambda\ddot\lambda
        \right)
    \right),
$$
$$
\lambda^2\frac{\partial^2y}{\partial t^2} =
\left(2\dot\lambda^2-\ddot\lambda \lambda\right)x\frac{Y}{X}
 +\frac{\dot\lambda x^3}{2}\lambda^3\frac{\partial}{\partial
t}\left[\frac{\lambda\ddot\lambda}{\lambda^2}\frac{3 A + 1 -\alpha + \frac{\alpha x}{
y}}{X}\right] $$
\begin{equation}
\label{31}
+ O\left(x^3\lambda\dot\lambda \frac{\partial}{\partial
t}\left(\lambda\ddot\lambda\right)\right)
\end{equation}
 where
\begin{multline}
\label{31_1}
\frac{\partial}{\partial t}\left[
\frac{\lambda\ddot\lambda}{\lambda^2}\frac{3A + 1 - \alpha + \frac{\alpha x}{ y}}{X}\right] \\
 =\frac{\dot\lambda\lambda\ddot\lambda}{\lambda^3}\biggl\{\frac{\alpha x^2}{y^2}\frac{Y^2}{X^3} - \frac{6\alpha x}{ y X^2} - \frac{3(1-\alpha)}{X}   -\frac{2 (2Z-1)\left(3A + 1 - \alpha\right)}{X^2}\biggr\}.
\end{multline}

In the domain $y \ge y_{cr}$ we have
$$
\lambda\frac{\partial y}{\partial t} = \frac{1}{Z}
\left\{
    \dot\lambda x Y + O
    \left(
        x^3 \lambda
        \frac{\partial}{\partial t} \left(
            \lambda \ddot\lambda
        \right)
    \right)
\right\},
$$
\begin{equation}
\label{31_2}
\lambda^2\frac{\partial^2 y}{\partial t^2} = -\left(2\dot\lambda^2 x -
\lambda \ddot\lambda x\right)\frac{Y}{Z}
 - \frac{\dot\lambda x^3}{2}\lambda^3\frac{\partial}{\partial t}\left[\frac{\lambda\ddot\lambda}{\lambda^2}
\frac{3A + 1 - \alpha -
\frac{\alpha x}{ y}}{Z}\right]
$$
$$
+
O\left(x^3\lambda\dot\lambda \frac{\partial}{\partial
t}\left(\lambda\ddot\lambda\right)\right),
\end{equation}
 where
\begin{equation}
\label{31_3}
\frac{\partial}{\partial t}\left[\frac{\lambda\ddot\lambda}{\lambda^2}
\frac{3A + 1 - \alpha -
\frac{\alpha x}{ y}}{Z}\right] $$$$=
\frac{\dot\lambda\lambda\ddot\lambda}{\lambda^3} \biggl \{\frac{\alpha x^2}{y^2}
\frac{Y^2}{Z^3} +
\frac{6\alpha x}{ y Z^2} -
\frac{3(1-\alpha)}{Z} - \frac{2 (2X-1)\left(3A + 1 -
\alpha\right)}{Z^2}\biggr\}.
\end{equation}

Now we present an explicit form of the function $\chi$ entering the definition of $\psi$, \eqref{29}, and introduced in \eqref{27'}. Due to Eqn \eqref{20} we have
\begin{equation}
\label{90} \chi = \begin{cases} -\lambda\ddot\lambda x^2
\left(A+1/2\right)X^{-1}, & x < x_{cr}  \\
\biggl[-\frac{2 y_{cr}}{ x} +
\lambda\ddot\lambda x^2\left(A+1/2\right)\biggr]
Z^{-1}, & x > x_{cr}.
\end{cases}
\end{equation}

Next, we give here an explicit expression for the expression $\partial\chi/\partial x - 2\chi/x$.
We compute
\begin{multline}
\label{30}
\frac{\partial\chi}{\partial x} - \frac{2\chi}{x} = -
\frac{\lambda\ddot\lambda x^2}{ X}
\biggl\{\frac{1-\alpha}{ x} + \frac{\alpha Y}{ y X}
 -
\frac{\alpha\lambda\ddot\lambda x^2}{2 y X} \left(A+1/2\right)\left(3-\frac{x}{
y}\frac{Y}{X}\right) \biggr\}, \text{}
\end{multline}
for  $x < x_{cr}$, and
\begin{multline}
\label{30_1}
\frac{\partial\chi}{\partial x}-\frac{2\chi}{x}=
\frac{6y_{cr}}{ x^2 Z} +
\frac{\lambda\ddot\lambda x^2}{ Z}
\left(\frac{1-\alpha}{ x} - \frac{\alpha}{y}\frac{Y}{Z}\right)\\ +
\frac{\alpha\lambda\ddot\lambda x}{2y Z^2}
\left(-2y_{cr}+\lambda\ddot\lambda x^3\left(A+1/2\right)\right)
\left(3+\frac{x}{y}\frac{Y}{Z}\right), 
\end{multline}
for $x > x_{cr}$.

%
Note, that the function $y$, Eqn \eqref{20}, is chosen so as to cancel the term
$-2\lambda \ddot\lambda x/(\lambda^2(1+y^2))$ arising from the last term in expression \eqref{29} (see the first term on the r.h.s. of \eqref{31} and the first term on the r.h.s. of \eqref{14}).
With the help of Eqns \eqref{29}, \eqref{30} and \eqref{31} we obtain following expression for the function
$\psi$ in the domain $x \le x_{cr}$:
\begin{equation}
\label{A3.1}
 \psi 
 =
\frac{x^2}{y^2}\biggl\{-\frac{8\lambda\ddot\lambda x}{(1+y^2)^2}
\frac{A+1/2}{X} - \frac{4\dot\lambda^2 x}{(1+y^2)^2 X}+\frac{2\dot\lambda^2\lambda\ddot\lambda x^4y\left(A - \alpha\right)}{(1 + y^2)^2 X^2} $$
$$ +
\frac{2\alpha\lambda\ddot\lambda x^2}{yX(1 + y^2)}\biggl[\frac{y}{x } - \frac{Y}{  X} +
\frac{\lambda\ddot\lambda x^2}{2 }
\frac{A + 1/2}{X}\left(3 - \frac{x}{ y}\frac{Y}{X}\right)\biggr] $$
$$- \frac{4(\lambda\ddot\lambda)^2 x^4y}{(1 + y^2)^2}
\frac{\left(A +
1/2\right)^2}{X^2} -
\frac{4\dot\lambda^2 y x^2(1 - Y^2)}{(1 + y^2)^2X^2} $$
$$- \frac{2x}{1 +y^2}\biggl[-\left(2\dot\lambda^2 -
\lambda\ddot\lambda \right)\frac{1-Y}{X}
+
\frac{\dot\lambda^2\lambda\ddot\lambda x^2}{2}\biggl(\frac{\alpha x^2}{y^2}\frac{Y^2}{X^3} - \\ \frac{6\alpha x}{ y\left(1 +
\frac{\alpha\lambda\ddot\lambda x^3}{2 y}\right)^2}  $$
$$ - \frac{3(1-\alpha)}{X}- \frac{2(2Z-1)\left(3A + 1 -
\alpha\right)}{X^2}\biggr)\biggr]\biggr\}.
\end{equation}

Now, using Eqns \eqref{90}, \eqref{29}, \eqref{30}, \eqref{31_2}, we find expression for function $\psi$ in the
domain $y \ge y_{cr}$. 
In fact, to obtain the equation on the parameter $\lambda$ we need to know only the part of $\psi$ in $\{y \ge y_{cr}\}$,
proportional to $\dot\lambda^2 
$. For this reason we write out only this part:
\begin{multline}
\label{A3.2}
\psi = 2\dot\lambda^2\frac{x^3}{y^4}\biggl\{\frac{2x}{y}\frac{Y^2}{Z^2} +
 \frac{2Y}{Z}
  +\frac{\lambda\ddot\lambda x^2}{2Z}\biggl(\frac{\alpha x^2}{y^2} \frac{Y^2}{Z^2} + \frac{6\alpha x}{ yZ}  - 3(1-\alpha)\\  -\frac{2 (2X-1)\left(3A +1 -
\alpha\right)}{Z}\biggr)\biggr\}
+ \mbox{ term proportional to }\ \lambda\ddot\lambda.
\end{multline}

Finally, we show Eqn \eqref{asbehavpsi}
which was stated in Section \ref{approximatesolution}. Indeed, the definitions of $Y$ and $Z$ and the second equation in \eqref{20} imply that for $ y \ge y_{cr}$
$$\frac{Y}{Z}=\frac{3\lambda\ddot\lambda x^2}{2\alpha}B^3,\ Y\sim \lambda\ddot\lambda x^2 B,\ y \sim \lambda\ddot\lambda x^3 B,$$
where $B:=\ln (\lambda\ddot\lambda x^2)$. Using these relations and the equation \eqref{A3.2},
we arrive at the desired relation \eqref{asbehavpsi}.

\bigskip

\textbf{Appendix} 4

\bigskip

In this appendix we compute the partial derivatives of energy $E=E(\alpha, \beta)$ w.r.to parameters $\alpha, \beta$. Using  expression \eqref{57}, we obtain
\begin{multline}
\label{58}
\frac{\partial E}{\partial \beta} = 4 \int\limits_0^\infty d\rho\frac{\rho}{(1+y^2)^2}\biggl\{\frac{\partial
y}{\partial t}\frac{\partial}{\partial \beta}\left(\frac{\partial y}{\partial t}\right) + \frac{\partial
y}{\partial\rho}\frac{\partial}{\partial \beta}\left(\frac{\partial y}{\partial\rho}\right) + \frac{y}{\rho^2}\frac{\partial y}{\partial \beta}  \\-
\frac{2y}{1+y^2}\frac{\partial y}{\partial \beta}\left(\left(\frac{\partial y}{\partial t}\right)^2 +
\left(\frac{\partial y}{\partial\rho}\right)^2 + \frac{y^2}{\rho^2}\right)\biggr\}
\end{multline}
and
\begin{multline}
\label{58_1}
\frac{\partial E}{\partial\alpha} = 4\int\limits_0^\infty d\rho \frac{\rho}{(1+y^2)^2}\biggl\{\frac{\partial
y}{\partial t}\frac{\partial}{\partial\alpha}\left(\frac{\partial y}{\partial t}\right) + \frac{\partial
y}{\partial \rho}\frac{\partial}{\partial\alpha}\left(\frac{\partial y}{\partial\rho}\right) +
\frac{y}{\rho^2}\frac{\partial y}{\partial\alpha}  \\- \frac{2y}{1+y^2} \left(\left(\frac{\partial
y}{\partial t}\right)^2 + \left(\frac{\partial y}{\partial\rho}\right)^2 + \frac{y^2}{\rho^2}\right)\biggr\}.
\end{multline}

Recall the notation
\begin{equation}
\label{XZ'}
X:=1 +\frac{\alpha\lambda\ddot\lambda x^3}{2y} \ \mbox{and}\ Z:=1 - \frac{\alpha\lambda\ddot\lambda x^3}{2 y}
\end{equation}
and let
\begin{equation}
\label{XZcr}
X_{cr}:=1 + \frac{\alpha\lambda\ddot\lambda x_{cr}^3}{2y_{cr}} \ \mbox{and}\ Z_{cr}:=1 - \frac{\alpha\lambda\ddot\lambda x_{cr}^3}{2 y_{cr}}.
\end{equation}
From Eqn \eqref{20} we find
\begin{align}
\label{A4.1}
\frac{\partial y}{\partial \beta} = \frac{\lambda\ddot\lambda x^3}{2 \beta}\frac{1}{X}, \qquad x < x_{cr}
\end{align}
$$
\frac{\partial y}{\partial \beta} = \frac{\lambda\ddot\lambda}{2 \beta} \frac{1}{Z}\left[\frac{2x_{cr}^3}{X_{cr}} - x^3\right], \qquad x > x_{cr}
$$
and
$$
\frac{\partial y_{cr}}{\partial \beta} = \frac{\lambda\ddot\lambda x_{cr}^3}{2 \beta} \frac{1}{X_{cr}}.
$$
Using Eqn \eqref{A4.1}, we obtain the time derivative of  $\frac{\partial y}{\partial \beta}$ in
the leading approximation in $1/\ln\left(\frac{1}{\lambda\ddot\lambda}\right)$:
\begin{equation}
\label{A4.2}
\frac{\partial}{\partial t}\left(\frac{\partial y}{\partial \beta}\right) =
-\frac{\dot\lambda\lambda\ddot\lambda x^3}{2\beta\lambda Z^2} \left[3 +
\frac{\alpha\lambda\ddot\lambda x^4}{2y^2} \frac{Y}{X}\right], \qquad x < x_{cr},
\end{equation}
\begin{multline}
\label{A4.3}
\frac{\partial}{\partial t}\left(\frac{\partial y}{\partial \beta}\right) =
-\frac{\dot\lambda\lambda\ddot\lambda}{2\beta\lambda Z}  \left\{ \frac{3 + \frac{\alpha\lambda\ddot\lambda x^4}{2y^2}\frac{Y}{Z}}{Z}  \left(\frac{2x_{cr}^3}{X_{cr}}
-x^3\right) - \frac{6x_{cr}^3}{X_{cr}}
\right\},\ 
x > x_{cr}.
\end{multline}
In a similar way we find derivative of $y$ w.r. to $\alpha$:
\begin{equation}
\label{A4.4}
\frac{\partial y}{\partial\alpha} = -\frac{\lambda\ddot\lambda x^3}{2}
\frac{\ln\left(\frac{ y}{x}\right)}{X},
\qquad x < x_{cr},
\end{equation}
$$
\frac{\partial y}{\partial\alpha} = -\frac{\lambda\ddot\lambda}{2} \frac{1}{Z} \left[-x^3\ln\left(\frac{ y}{x}\right) +
2x_{cr}^3 \frac{\ln\left(\frac{ y_{cr}}{x_{cr}}\right)}{X_{cr}}\right], \qquad x > x_{cr},
$$
$$
\frac{\partial y_{cr}}{\partial\alpha} = -\frac{\lambda\ddot\lambda x_{cr}^3}{2}
\frac{\ln\left(\frac{ y_{cr}}{x_{cr}}\right)}{X_{cr}}.
$$
Taking the time derivative of Eqn \eqref{A4.4}, we obtain
\begin{multline}
\label{A4.5}
\frac{\partial}{\partial t}\left(\frac{\partial y}{\partial\alpha}\right) =
-\frac{\dot\lambda\lambda\ddot\lambda x^3}{2\lambda X^2}
\biggl\{\frac{\lambda\ddot\lambda x^3}{ y}
\left(\ln\left(\frac{\sqrt{\lambda\ddot\lambda}}{\beta}y^\alpha x^{1-\alpha}\right) +1/2\right) - \\
\ln\left(\frac{ y}{x}\right) \left(3 + \frac{\alpha\lambda\ddot\lambda x^4}{2y^2}
\frac{Y}{X}\right)\biggr\}, \qquad x < x_{cr},
\end{multline}
\begin{multline}
\label{A4.6}
\frac{\partial}{\partial t}\left(\frac{\partial y}{\partial\alpha}\right) =
-\frac{\dot\lambda\lambda\ddot\lambda}{2\lambda Z} \biggl\{ - \frac{3 + \frac{\alpha  \lambda\ddot\lambda x^4}{2 y} \frac{Y}{Z}}{Z}
\biggl[2x_{cr}^3 \frac{\ln\left(\frac{ y_{cr}}{x_{cr}}\right)}{X_{cr}} - x^3\ln\left(\frac{
y}{x}\right)\biggr]  + \frac{6x_{cr}^3\ln\left(\frac{ y_{cr}}{x_{cr}}\right)}{X_{cr}}\\ - \frac{x^3}{y Z}  \biggl(2y_{cr} -
\lambda\ddot\lambda x^3 \left(\ln\left(\frac{\sqrt{\lambda\ddot\lambda}}{\beta}
y^\alpha x^{1-\alpha}\right) +1/2\right)\biggr)\biggr\}, \qquad x > x_{cr}.
\end{multline}

Note that the main contribution to the partial derivatives \eqref{58} and \eqref{58_1} comes from the domain $x\sim x_{cr}$. Both
derivatives are sums of terms proportional to $\dot\lambda^2$ and  to
$\lambda\ddot\lambda$. The coefficients for these terms are of order of $1$. As result, since we assumed that $|\lambda\ddot\lambda| \ll \dot\lambda^2$, we have to
find in the expressions for \eqref{58} and \eqref{58_1} only the terms proportional to $\dot\lambda^2$.

Using Eqns \eqref{A4.1} - \eqref{A4.6} we can write the r.h.s. of \eqref{58} and \eqref{58_1} in a more explicit form
\begin{multline}
\label{A4.7}
\frac{1}{4}\frac{\partial E}{\partial \beta} = \dot\lambda^2 \int\limits_0^{x_{cr}}
\frac{dx x^5\lambda\ddot\lambda}{2\beta y^4} \frac{Y}{\left(1 +
\frac{\alpha\lambda\ddot\lambda x^3}{2 y}\right)^3} \biggl[ 3 +
\frac{\alpha\lambda\ddot\lambda x^4}{2 y^2} \frac{Y}{1+\frac{
\alpha\lambda\ddot\lambda x^3}{2 y}} - \frac{2x Y}{ y}\biggr]  \\
- \dot\lambda^2
\int\limits_{x_{cr}}^\infty \frac{dx x^2\lambda\ddot\lambda}{2\beta y^4} \frac{Y}{\left(1 -
\frac{\alpha\lambda\ddot\lambda x^3}{2 y}\right)^2} \biggl\{\biggl[\frac{1}{1 -
\frac{\alpha\lambda\ddot\lambda x^3}{2 y}} \biggl(3 +
\frac{\alpha\lambda\ddot\lambda x^4}{2 y^2} \frac{Y}{1 -
\frac{\alpha\lambda\ddot\lambda x^3}{2 y}}\biggl) \\ \times \left(\frac{2x_{cr}^3}{1 +
\frac{\alpha\lambda\ddot\lambda x_{cr}^3}{2 y_{cr}}} - x^3\right)
- \frac{6x_{cr}^3}{1 +
\frac{\alpha\lambda\ddot\lambda x_{cr}^3}{2 y_{cr}}}\biggr] + \frac{2x}{ y} \frac{Y}{1 -
\frac{\alpha\lambda\ddot\lambda x^3}{2 y}} \left(\frac{2x_{cr}^3}{1 +
\frac{\alpha\lambda\ddot\lambda x_{cr}^3}{2 y_{cr}}} -x^3\right)\biggr\},
\end{multline}
\begin{equation}
\label{A4.8}
\frac{1}{4}\frac{\partial E}{\partial\alpha} = \dot\lambda^2 \int\limits_0^{x_{cr}}
\frac{dx x^5\lambda\ddot\lambda Y}{2 y^4\left(1 +
\frac{\alpha\lambda\ddot\lambda x^3}{2 y}\right)^3}
\biggl\{\frac{\lambda\ddot\lambda x^3}{ y}\left(\ln\left(\frac{\sqrt{\lambda\ddot\lambda}}{\beta}
y^\alpha x^{1-\alpha}\right) + 1/2\right) $$
$$-
\ln\left(\frac{ y}{x}\right) \left(3 +
\frac{\alpha\lambda\ddot\lambda x^4}{2 y^2} \frac{Y}{1 +
\frac{\alpha\lambda\ddot\lambda x^3}{2 y}}\right) + \frac{2x}{ y} \ln\left(\frac{
y}{x}\right) Y\biggr\} + \dot\lambda^2 \int\limits_{x_{cr}}^\infty \frac{dx x^2\lambda\ddot\lambda
Y}{2 y^4 \left(1 - \frac{\alpha\lambda\ddot\lambda x^3}{2 y}\right)^2} $$
$$\times
\biggl\{\frac{1}{1 - \frac{\alpha\lambda\ddot\lambda x^3}{2 y}} \left(3 +
\frac{\alpha\lambda\ddot\lambda x^4}{2 y^2} \frac{Y}{1 -
\frac{\alpha\lambda\ddot\lambda x^3}{2 y}}\right) \biggl(\frac{2x_{cr}^3\ln\left(\frac{y_{cr}}{x_{cr}}\right)}{1 + \frac{\alpha\lambda\ddot\lambda x_{cr}^3}{2 y_{cr}}} - x^3
\ln\left(\frac{ y}{x}\right)\biggr) $$
$$- \frac{6x_{cr}^3 \ln\left(\frac{y_{cr}}{x_{cr}}\right)}{1 + \frac{\alpha\lambda\ddot\lambda x_{cr}^3}{2 y_{cr}}} +
\frac{x^3}{y\left(1 - \frac{\alpha\lambda\ddot\lambda x^3}{2 y}\right)} \left(2y_{cr} -
\lambda\ddot\lambda x^3 \left(\ln\left(\frac{\sqrt{\lambda\ddot\lambda}}{\beta} y^\alpha
x^{1-\alpha}\right) + 1/2\right)\right) $$
$$+ \frac{2x}{ y} \frac{Y}{1
-\frac{\alpha\lambda\ddot\lambda x^3}{2 y}} \left(\frac{2x_{cr}^3\ln\left(\frac{
y_{cr}}{x_{cr}}\right)}{1 + \frac{\alpha\lambda\ddot\lambda x_{cr}^3}{2 y_{cr}}} - x^3
\ln\left(\frac{ y}{x}\right)\right) \biggr\}.
\end{equation}

Let  $\alpha = \alpha( \beta)$ be a  solution of Eqn \eqref{54}. We find numerically (see Appendix 5)
that $\beta$ changes on the interval $(0, \beta_0]$, where $\beta_0=1.0405$ (the corresponding value of $\alpha$ is $\alpha_0 = 0.65436.$) Using expressions for $\frac{\partial E}{\partial\alpha} $ and $\frac{\partial E}{\partial\beta}$, derived above, we show numerically that  the function
\begin{equation}
\label{61_1}
\Phi := \frac{\partial E}{\partial \beta} + \frac{\partial E}{\partial\alpha}\frac{\partial\alpha}{\partial \beta}
\end{equation}
 is negative for $\beta =\beta_0$ and for $\beta \rightarrow 0$, with $E( \alpha, \beta)$  having absolute minimum at $\beta =\beta_0$.

 \bigskip

\textbf{Appendix} 5

\bigskip

 Numerical calculations with help of Eqns \eqref{55} and \eqref{55_1} show that there is a point $(\alpha_0,\beta_0)$,
 $$\beta_0=1.0405, \alpha_0=.6543626,$$
  so that the equation  \eqref{54}, $ I(\alpha,\beta)=0 $, has no solution for $\alpha<\alpha_0$ and  for $\beta > \beta_0$. Moreover, the solution of the equation  $ I(\alpha,\beta)=0 $ for  $\beta$ determines a double-valued function   $ \beta=\beta(\alpha)$, whose branches coalesce at $\alpha=\alpha_0$ and have different derivatives there (see Eqns \eqref{A5.1} and \eqref{Deltaexp} below).
 Moreover,  $ I(\alpha,\beta)=0 $ has the unique solution $\beta_0$ at $\alpha=\alpha_0$.  Hence the solution of the equation  $ I(\alpha,\beta)=0 $ for  $\alpha$ also leads to a double-valued function $\alpha=\alpha (\beta)$.

 Numerical calculations give the following expansions for the lowest branch,
 \begin{equation}
  \label{A5.1}
 \beta=\beta_0-
 \beta_1(\alpha-\alpha_0)-\beta_2(\alpha-\alpha_0)^2,
 \end{equation}
  $\alpha>\alpha_0$, and and for  the distance, $\Delta$, between the branches along the $\alpha$-axis,
 \begin{equation} \label{Deltaexp}
 \Delta=\gamma_1(\beta_0-\beta)-\gamma_2(\beta_0-\beta)^2,
 \end{equation}
where
\begin{align}
\label{A5.2}
 \beta_1=2.54732,\ \beta_2=13.8297, \\ \gamma_1=.08029,\ \gamma_2=.42736. 
\end{align}
(Solving \eqref{A5.1} for $\alpha$ gives the lower branch of the function  $\alpha=\alpha (\beta)$. Adding \eqref{Deltaexp} to this solution gives the upper branch of  $\alpha=\alpha (\beta)$.)

To find the second "end" point on the $\alpha-$interval we check the point  $\alpha=1$
where the dependence of $y $ on $x$ in Eqn \eqref{20} can be found in an explicit form. To do this we note that \eqref{23} and \eqref{24} with $\alpha =1$ imply that
 \begin{equation}
 \label{A5.3}
 \gamma= { \sqrt{\lambda\ddot\lambda} y_{cr} \over \beta },\
    \beta^2=\frac{8}{27 \gamma^2\ln \gamma}\ \mbox{and}\
 \lambda\ddot\lambda {y_{cr}}^2=\frac{8}{27\ln \gamma}. 
 \end{equation}
We also have $ y_{cr}/x_{cr}=2/3$.  For $\alpha=1$ solvability condition of Eqn \eqref{20} is
 \begin{equation}
  \label{A5.5}
 \gamma>e^{1/2}.
 \end{equation}
  Indeed, set
  \begin{equation}
  \label{A_5}
  y=y_{cr} z,\;\;z=1+\delta,\;\; \frac{x}{y_{cr}}= \frac{3}{2}+\tau
  \end{equation}
  In the range $0<\delta\ll 1$ we have
  \begin{equation}
  \label{A_6}
  \frac{2}{3} \tau^2=\delta \left(1- \frac{1}{2\ln{\gamma}}\right)
  \end{equation}
  From this equation we see that $\beta$ should satisfy the inequality given in Eqn \eqref{A5.5}.

 Now we set $y=y_{cr}z$. For $z<1$ we obtain from the first equation in \eqref{20}, with $\alpha=1$, and from \eqref{A5.3} the following cubic equation for the ratio  $\frac{x}{ y_{cr}}$
 \begin{equation}
 \label{A5.7}
 \frac{4}{27}\biggl(\frac{x}{ y_{cr}}\biggr)^3\frac{\ln(\gamma z)}{\ln \gamma}-\frac{x}{ y_{cr}}+z=0.
 \end{equation}
 Solution of Eqn \eqref{A5.7} in the range $z<1$ is
 \begin{equation}
 \label{A5.8}
\frac{x}{ y_{cr}}= 3\biggl(
    \frac
        {\ln \gamma}
        {\ln (1/(\gamma z))}
\biggr)^{1/2} \sinh\biggl[
\frac{1}{3}
\ln\biggl(z
    \sqrt{\frac
        {\ln (1/(\gamma z))}
        {\ln \gamma}
    }+
    \sqrt{1+z^2\frac
        {\ln (1/(\gamma z))}
        {\ln \gamma}
    }
\biggr)\biggr]
\end{equation}
 for $\gamma z<1$ and
 \begin{equation}
 \label{A5.9}
 \frac{x}{ y_{cr}}=3\sqrt{\frac{\ln \gamma}{\ln (\gamma z)}}\sin\biggl[\frac{1}{3}\arctan\frac{z\sqrt{\ln (\gamma z)/\ln \gamma}}
{\sqrt{1-z^2\ln (\gamma z)/\ln \gamma}}\biggr]
 \end{equation}
 for $\gamma z>1$.

  In the range $z \ge 1$ the ratio $x/ y_{cr}$ solves the following cubic equation (see the second equation in \eqref{20})
  \begin{equation}
  \label{A5.10}
  \frac{4}{27}\frac{\ln(\gamma z)}{\ln \gamma}\biggl(\frac{x}{ y_{cr}}\biggr)^3-\frac{x}{ y_{cr}}+2-z=0.
  \end{equation}
  Let $z_0$ be solution of equation
  \begin{equation}
  \label{A5.13}
  1=(z_0-2)^2\frac{\ln(\gamma z_0)}{\ln \gamma}.
  \end{equation}
We split the semi-interval $z>1$ into two sub-intervals. In the interval $ 1<z<z_0$ we have
  \begin{equation}
   \label{A5.11}
  \frac{x}{ y_{cr}}=3\sqrt{\frac{\ln \gamma}{\ln(\gamma z)}}\sin\phi,
  \end{equation}
  where
  \begin{multline}
    \label{A5.12}
  \phi=\frac{\pi}{6}+\frac{1}{3}\arctan\biggl(
    \frac
        {\sqrt{1-(2-z)^2\ln (\gamma z)/\ln \gamma}}
        {(2-z)\sqrt{\ln (\gamma z)/\ln \gamma}}
  \biggr), 1<z<2, \\
  \phi=\frac{\pi}{3} +\frac{1}{3}\arctan\biggl(
    \frac
        {(z-2)\sqrt{\ln (\gamma z)/\ln \gamma}}
        {\sqrt{1-(z-2)^2\ln (\gamma z)/\ln \gamma}}
  \biggr), 2<z<z_0.
  \end{multline}
  In the range $z>z_0$ we have
  \begin{equation}
    \label{A5.14}
  \frac{x}{ y_{cr}}=\frac{3}{2}(Q^{1/3}+\frac{\ln \gamma}{\ln(\gamma z)}Q^{-1/3})
\end{equation}
where
\begin{equation}
\label{A5.15}
Q=(z-2)\frac{\ln \gamma}{\ln(\gamma z)}+\sqrt{\biggl((z-2)\frac{\ln \gamma}{\ \ln(\gamma z)}\biggr)^2-\biggl(\frac{\ln \gamma}{\ln(\gamma z)}\biggr)^3}
\end{equation}
 Using  Eqns \eqref{A5.8}-\eqref{A5.15} we obtain with help of numerical  calculations, that Eqn \eqref{54}
 at $\alpha=1$ has solution only as $\beta$ goes to zero. This means that  $\alpha=1$ is the second end point of the $\alpha-$interval.


\begin{thebibliography}{99}

\bibitem{At} m.F. Atiah, Geometry of Yang-Mills Fields. Accademia Nazionale Dei Lincei Scuola Normale Superiore, 1979.


\bibitem{Bi} P. Bizo\'n, Formation of singularities in Yang-Mills equations,
Acta Phys. Polonica B. 33, 1893 (2002).

\bibitem{BiPC} P. Bizo\'n, Private communication.
\bibitem{BCT} P. Bizo\'n, T. Chmaj, Z. Tabor, Formation of singularities for equivariant (2+1)-dimensional wave maps into the 2-sphere. Nonlinearity 14 (2001), no. 5, 1041-1053.

\bibitem{BOS} P. Bizo\'n, Yu. N. Ovchinnikov, I. M. Sigal, Collapse of an instanton. Nonlinearity {\bf 17} (2004),
no. 4, 1179-1191.
\bibitem{BT} P. Bizo\'n and Z. Tabor, On blowup of Yang-Mills fields,
Phys. Rev. D64, 121701 (2001).
\bibitem{B} E. B. Bogomolnyi, The stability of classical solutions. Soviet J. Nuclear Phys. 24 (1976), no.
4, 449-454 (Russian).
\bibitem{CST} T. Cazenave, J. Shatah, S. Tahvildar-Zadeh, Harmonic maps of the hyperbolic space and
development of singularities in wave maps and Yang-Mills fields. Ann. I.H.P., section A 68
(1998), no. 3, 315-349.
\bibitem{CDY} K.C. Chang, W. Y. Ding, R. Ye, Finite-time blow-up of the heat flow of harmonic maps from
surfaces. J. Differential Geom. 36 (1992), no. 2, 507-515.
\bibitem{CS} Y. M. Chen, M. Struwe, Existence and partial regularity results for the heat flow for harmonic
maps. Math. Z. 201 (1989), no. 1, 83-103.
\bibitem{CTZ1} D. Christodoulou, S. A. Tahvildar-Zadeh, On the regularity of spherically symmetric wave
maps. Comm. Pure Appl. Math. 46 (1993), no. 7, 1041-1091.
\bibitem{CTZ2} D. Christodoulou, S. A. Tahvildar-Zadeh, On the asymptotic behavior of spherically symmetric
wave maps. Duke Math. J. 71 (1993), no. 1, 31-69.
\bibitem{9} R. C\^ote, Instability of non-constant harmonic maps for the (1+2)-dimensional equivariant
wave map system. Int. Math. Res. Not. 57 (2005) 3525-3549.

\bibitem{GS} J.F. Grotowski, J. Shatah, Geometric evolution equations in critical dimensions. preprint.

\bibitem{Handbook} "Handbook of Mathematical functions" Edited by M.
Abramowitz and I.A. Stegun, National Burean of Standarts, Applied
Mathematical Series. 55 (1964).
\bibitem{IL}  J. Isenberg and S.L. Liebling, Singularity formation in 2+1 wave maps. J. Math. Phys. 43
(2002), 678-683.
\bibitem{KlMach} S. Kleinerman and M. Machedon, Smoothing estimates for null forms and applications. Internat. Math. res. Notes, 13:655-677, 2001.
\bibitem{KlSel1} S. Kleinerman and S. Selberg, Remark on optimal regularity for equations of wave maps type. Comm. Partial Differential Equations, 22 (2-5): 901-918, 1997.
\bibitem{KlSel2} S. Kleinerman and S. Selberg, Bilinear estimates and applications to nonlinear wave equations. Comm. Contemp. Math.,4 (2):223-295, 2002.
\bibitem{Kov} M. Kovalyov, Long-time behaviour of solutions of a system of nonlinear wave equations. Comm. Partial Differential Equations, 12 (5): 471-501, 1987.
\bibitem{Kr1} J. Krieger, Global regularity of wave maps from R2+1 to H2. Small energy. Comm. Math.
Phys. 250 (2004), no. 3, 507-580.
\bibitem{Kr2} J. Krieger, Stability of spherically symmetric wave maps. to appear.
\bibitem{KST1} J. Krieger, W. Schlag and D. Tataru, Renormalization and blow up for charge one equivariant critical wave maps, Invent. Math. 171 (2008), no. 3, 543-615.
\bibitem{KST2} J. Krieger, W. Schlag and D. Tataru, Renormalization and blow up for critical Yang-Mills problem, e-print, arXiv:0809.211, 2008.
\bibitem{LPZ} R. A. Leese, M. Peyrard, W. J. Zakrzewski, Soliton stability in the O(3) $\sigma$-model in (2 + 1)
dimensions. Nonlinearity 3 (1990), no. 2, 387-412.
\bibitem{LS} J.M. Linhart, L. Sadun, Fast and slow blowup in the S2 $\sigma$-model and the (4+1)-dimensional
Yang-Mills model. Nonlinearity 15 (2002), 219-238.
\bibitem{MS}  N. Manton, P. Sutcliffe, Topological Solitons. Cambridge Monographs on Mathematical
Physics. Cambridge University Press, Cambridge, 2004.

\bibitem{PZ} B. Piette, W. J. Zakrzewski, Shrinking of solitons in the (2+1)-dimensional S2 sigma model.
Nonlinearity 9 (1996), no. 4, 897-910.
\bibitem{RRS} P.Rafael, I.Rodnianski and J.Sternbenz, In preparation.
\bibitem{R} Rendall A.D., in "Applications of the theory of evolution
equations to general relativity", Proc. GR16 ed. N.T. Bishop and
S.D. Makaraj (Singapure: World Scientific).
\bibitem{RS} I.Rodnianski and J.Sternbenz, On the formation of singularities in the critical $O(3)\ \sigma$-model, e-print, arXiv:math/0605023, 2006.
\bibitem{Schw} A. S. Schwarz, Quantum Field Theory and Topology. Springer-Verlag, 1993.
\bibitem{Sh} J. Shatah, Weak solutions and development of singularities in the SU(2) $\sigma$-model. Comm.
Pure Appl. Math. 41 (1988), 459-469
\bibitem{ShTZ1} J. Shatah, S. A. Tahvildar-Zadeh, Regularity of harmonic maps from the Minkowski space
into rotationally symmetric manifolds. Comm. Pure Appl. Math. 45 (1992), no. 8, 947-971.
\bibitem{ShTZ2} J. Shatah, S. A. Tahvildar-Zadeh, On the Cauchy problem for equivariant wave maps. Comm.
Pure Appl. Math. 47 (1994), no. 5, 719-754.
\bibitem{Sid} T. Sideris, Global existence of harmonic maps in Minkowski space, Comm. Pure Appl. Math., 42 (1): 1-13, 1989.
\bibitem{Str1} M. Struwe, Radially symmetric wave maps from (1 + 2)-dimensional Minkowski space to the
sphere. Math. Z. 242 (2002), no. 3, 407-414.
\bibitem{Str2} M. Struwe, Radially symmetric wave maps from (1 + 2)-dimensional Minkowski space to
general targets. Calc. Var. Partial Differential Equations 16 (2003), no. 4, 431-437.
\bibitem{Str3} M. Struwe, Equivariant wave maps in two space dimensions. Comm. Pure Appl. Math. 56
(2003), 815-823.
\bibitem{Tao1} T. Tao, Global regularity of wave maps. I. Small critical Sobolev norms in high dimensions. Internat. Math. Res. Notes, (6): 299-328, 2001.
\bibitem{Tao2} T. Tao, Global regularity of wave maps. II. Small energy in two dimensions. Comm. Math.
Phys. 224 (2001), no. 2, 443-544.
\bibitem{Tat1} D. Tataru, Local and global results for wave maps I.  Comm. Partial Differential Equations, 23 (9-10): 1781-1793, 1998.
\bibitem{Tat2} D. Tataru, On global existence and scattering for the wave maps equation. Amer. J. Math.
123 (2001), no. 1, 37-77.
\bibitem{Tat3} D. Tataru, The wave maps equation. Bull. AMS 41 (2): 185-204, 2004.
\bibitem{Tat4} D. Tataru, Rough solutions for the wave maps equation. Amer. J. Math. 127 (2005), no. 2,
293-377.
\bibitem{W} R.S. Ward, Slowly-moving lumps in the $\mathbb{C}\mathbb{P}^1$ model in $(2+1)$ dimensions. Physics letters 158B(5):424-428, 1985.





\end{thebibliography}
\end{document}